\documentclass[reqno,12pt]{amsart}
\usepackage{amssymb}
\usepackage{amsmath}
\usepackage{bbm}
\usepackage{mathrsfs}
\def\newpic#1{%
   \def\emline##1##2##3##4##5##6{%
      \put(##1,##2){\special{em:point #1##3}}%
      \put(##4,##5){\special{em:point #1##6}}%
      \special{em:line #1##3,#1##6}}}
\newpic{}

\makeatletter
\renewcommand{\@biblabel}[1]{#1.}
\makeatother

\theoremstyle{plain}
\newtheorem{prp}{Proposition}
\newtheorem{lemma}[prp]{Lemma}
\newtheorem{thm}[prp]{Theorem}
\newtheorem{cor}[prp]{Corollary}
\theoremstyle{remark}
\newtheorem{remark}{Remark}
\numberwithin{prp}{section}

\newcommand{\rf}[1][]{\textup{\eqref{#1}}}
\newcommand{\ov}{\overline}
\newcommand{\plus}[1]{\widetilde{#1}}

\newcommand{\Rda}[2][]{\hat{R}^{#1}_{#2}}
\newcommand{\Rch}[2][]{\check{R}^{#1}_{#2}}
\newcommand{\Rre}[2][]{\acute{R}^{#1}_{#2}}
\newcommand{\Rli}[2][]{\grave{R}^{#1}_{#2}}
\newcommand{\Rlim}[2][]{\grave{R}^{-#1}_{#2}}
\newcommand{\Rdam}[1][]{\hat{R}^{-1}_{#1}}

\newcommand{\Rchm}[1][]{\check{R}^{-1}_{#1}}
\newcommand{\pda}[1][\tau]{\hat{P^{#1}}}
\newcommand{\pch}[1][\tau]{\check{P^{#1}}}

\DeclareMathOperator{\id}{id}
\DeclareMathOperator{\End}{End}
\DeclareMathOperator{\Mor}{Mor}
\DeclareMathOperator{\im}{im}
\newcommand{\adr}{\mathrm{ad}_\rr}
\newcommand{\Adr}{\mathrm{Ad}_\rr}
\DeclareMathOperator{\rank}{rk}
\DeclareMathOperator{\tr}{tr}
\newcommand{\trq}{\tr^1_q}

\newcommand{\xx}{x}
\newcommand{\cv}{c}
\newcommand{\dd}{\mathrm{d}}
\newcommand{\ot}{\otimes}
\newcommand{\ott}{{\otimes}}
\newcommand{\ota}{{\otimes}_{\!{\scriptscriptstyle\mathcal{A}}}}
\newcommand{\ww}{{\mbox{$\scriptscriptstyle{\mathrm W}$}}}
\renewcommand{\SS}{\mathscr{S}}
\renewcommand{\AA}{\mathcal{A}}
\newcommand{\BB}{\mathcal{B}}
\newcommand{\RR}{\mathcal{R}}
\newcommand{\XX}{\mathcal{X}}
\newcommand{\C}{\mathbbm{C}}
\newcommand{\NO}{\mathbbm{N}_0}
\newcommand{\slqn}[1][N]{\mathrm{SL}_q(#1)}
\newcommand{\glqn}[1][N]{\mathrm{GL}_q(#1)}
\newcommand{\oqn}[1][N]{\mathrm{O}_q(#1)}
\newcommand{\spqn}[1][N]{\mathrm{Sp}_q(#1)}
\newcommand{\ooq} [1][N]{\mathcal{O}(\mathrm{O}_q(#1))}
\newcommand{\ospq}[1][N]{\mathcal{O}(\mathrm{Sp}_q(#1))}
\newcommand{\nn}{\nonumber}
\newcommand{\qm}[1][1]{q^{-{#1}}}
\newcommand{\Q}{Q}
\newcommand{\QM}{\hat{q}}
\newcommand{\RM}{\hat{\rt}}
\newcommand{\lam}{{\lambda}}
\newcommand{\eps}{\epsilon}
\newcommand{\ve} {\varepsilon}
\newcommand{\vem}{\varepsilon_-}
\newcommand{\vepm}{\varepsilon_\pm}
\newcommand{\thq}{\ov{\theta}}
\newcommand{\etat}[1][\tau]{\eta^{#1}}
\newcommand{\alpt}[1][\tau]{\alpha_{#1}}
\newcommand{\lamt}[1][\tau]{\lambda_{#1}}
\newcommand{\lamn}[1][\nu]{\lambda_{#1}}
\newcommand{\tnn}[1][\nu]{t_{#1}}
\newcommand{\ctn}[1][\tau\nu]{c_{#1}}
\newcommand{\vr}{\varrho}
\newcommand{\om}{\omega}
\newcommand{\fettc}{{\scriptstyle{\mathrm{c}}}}
\newcommand{\rac}{{\,\triangleleft\,}}
\newcommand{\fettt}{{\scriptscriptstyle{\top}}}

\newcommand{\uhr}{{\upharpoonright}}
\newcommand{\ii}{{\scriptscriptstyle{\mathrm{I}}}}
\renewcommand{\ll}{{\scriptscriptstyle{\mathrm{L}}}}
\newcommand{\rr}  {{\scriptscriptstyle{\mathrm{R}}}}
\newcommand{\rt}  {r}
\newcommand{\rtm}  {r^{-1}}
\newcommand{\dl}{{\Delta_\ll}}
\newcommand{\dr}{{\Delta_\rr}}
\newcommand{\qtwo}{[2]_q}
\newcommand{\Gamm}{\varGamma}
\newcommand{\Lam}[1][]{\varLambda^{#1}}
\newcommand{\gpm}{\Gamm_\pm}
\newcommand{\gl}{{\Gamm_\ll}}
\newcommand{\gr}{{\Gamm_\rr}}
\newcommand{\gi}{{\Gamm_\ii}}
 
\newcommand{\sig}{\sigma}
\newcommand{\gdw}[1][]{{_\ww{\Gamm^{\land #1}}}}
\newcommand{\gds}[1][]{{_{\mathrm s}{\Gamm^{\land #1}}}}
\newcommand{\gdu}[1][]{{_{\mathrm u}{\Gamm^{\land #1}}}}

\newcommand{\Js}{{_{\mathrm s} J}}
\newcommand{\Ju}{{_{\mathrm u} J}}

\newcommand{\gdsl}[1][]{{_{\mathrm s}{\Gamm^{\land #1}_\ll}}}
\newcommand{\gdsi}[1][]{{_{\mathrm{s}}{\Gamm^{\land #1}_\ii}}}
\newcommand{\gdul}[1][]{{_{\mathrm u}{\Gamm^{\land #1}_\ll}}}

\newcommand{\gdui}[1][]{{_{\mathrm u}{\Gamm^{\land {#1}}_\ii}}}

\newcommand{\gten} [1][] {\Gamm^{\ot #1}}

\newcommand{\gteni}[1][]{{\Gamm^{\ot #1}_\ii}}
\newcommand{\pair}{{\langle\cdot,\cdot\rangle}}
\newcommand{\ssp}{\discretionary{}{}{\,}}

\begin{document}
\author{Axel Sch\" uler}

\thanks{Supported by the Deutsche Forschungsgemeinschaft, e-mail: schueler@mathematik.uni-leipzig.de}

\address{Axel Sch\" uler, Department of Mathematics, 
University of Leipzig, Augustusplatz 10, 04109 Leipzig,
Germany}

\title[Exterior algebras for quantum groups]{Two exterior algebras for orthogonal and symplectic quantum groups}

\begin{abstract}
Let $\Gamm$ be one of the $N^2$-dimensional bicovariant first order
differential calculi on the quantum groups $\oqn$ or $\spqn$, where $q$ is not
a root of unity. 
\\
We show that the second antisymmetrizer exterior algebra $\gds$ is the
quotient of the universal exterior algebra $\gdu$ by the principal ideal
generated by $\theta{\land}\theta$. Here $\theta$ denotes the unique up to
scalars bi-invariant 1-form. Moreover $\theta{\land}\theta$ is central in
$\gdu$ and $\gdu$ is an inner differential calculus.
\end{abstract}

\maketitle


\section{Introduction}
More than a decade ago Woronowicz provided a general framework for covariant
differential calculus over arbitrary Hopf algebras, \cite{a-Woro2}. In the meantime there
exists a well developed theory of covariant differential calculus on Hopf
algebras, \cite[Chapter\,14]{b-KS}. In his paper Woronowicz also introduced the concept of
higher order forms which is based on a braiding
$\sig\colon\Gamm\ota\Gamm\to\Gamm\ota\Gamm$. The braiding $\sig$ naturally
generalizes the classical flip automorphism. It turns out that Woronowicz'
external algebra $\gdw$ is not simply a bicovariant bimodule but
a differential graded Hopf algebra \cite{a-Drabant1}, \cite[Theorem\,14.17]{b-KS}. However there
are two other concepts of exterior algebras which are also
differential graded Hopf algebras, \cite{a-Brz1,a-LyuSud1},
\cite[Theorem\,14.18]{b-KS}. The  ``second antisymmetrizer" exterior algebra
$\gds$ is also constructed using the braiding; but it involves only the
antisymmetrizer $I-\sig$ of second degree while Woronowicz' construction uses
antisymmetrizers of all degrees. The universal exterior algebra $\gdu$ can be
characterized by the following universal property. Each (left-covariant)
differential calculus which contains a given first order differential calculus
$\Gamm$ as its first order part is a quotient of $\gdu$. It seems  natural to
ask for the relation between these three concepts. For the quantum groups
$\glqn$ and $\slqn$ and their standard bicovariant first order differential calculi (abbreviated FODC) this problem was completely solved in \cite{a-Sch1}.
\\
In this paper we consider the quantum groups $\oqn$ and $\spqn$ together with
their standard bicovariant FODC.  The main result
is stated in Theorem\,\ref{t-1}. Suppose that $q$ is not a root of unity and
let $\theta$ be the unique up to scalars bi-invariant 1-form of $\Gamm$. Then
$\gdu{/}(\theta^2)$ and $\gds$ are isomorphic differential graded Hopf
algebras. Further, $\theta^2$ is central in $\gdu$ and $\gdu$ is an inner
differential calculus i.\,e. $\dd
\rho=\theta{\land}\rho-(-1)^n\rho{\land}\theta$ for $\rho\in\gdu[n]$. It is
somehow astonishing that the left-invariant parts of $\gdu[2]$ and $\gds[2]$
differ only  by the single element $\theta^2$.
\\
This paper is organized as follows. Section\,\ref{prelim} contains general
notions and facts about bicovariant bimodules and bicovariant differential
calculi over Hopf algebras. In Section\,\ref{osp} we recall necessary facts
about morphisms of corepresentations for orthogonal and symplectic quantum
groups. We give a brief introduction into the graphical calculus with
morphisms. The construction of bicovariant FODC on orthogonal and symplectic
quantum groups is reviewed. The main result is stated in Theorem\,\ref{t-1}.
 In Section\,\ref{proof}
a very useful criterion for the size
of the space of left-invariant 2-forms of $\gdu$ in terms of the quantum Lie algebra 
is given. This criterion applies to arbitrary left-covariant differential calculi.
We show that $\Gamm\ota\Gamm$ is the direct sum of 9 bicovariant subbimodules. Every  bicovariant subbimodule of $\Gamm\ota\Gamm$ which contains $\theta\ota\theta$
already contains the kernel of $I-\sig$.  Section\,\ref{uni} exclusively deals with the universal differential
calculus. The  outcome of the very technical calculations is that
$\theta{\land}\theta$ is non-zero and the unique up to scalars bi-invariant
2-form in $\gdu$.

We close the introduction by  fixing assumptions and notations that are used
in the sequel.  All vector spaces,
algebras, bialgebras, etc.\ are meant to be $\C$-vector spaces, unital $\C$-algebras,
$\C$-bialgebras etc. The linear span of a set $\{a_i| i\in K\}$ is
denoted by $\langle a_i| i\in K\rangle$.
$\AA$ always denotes  a  Hopf algebra. We write $\AA^\circ$ for the dual Hopf algebra. 
All modules, comodules, and bimodules are assumed to be $\AA$-modules, 
$\AA$-comodules, and $\AA$-bimodules if nothing else is specified. 
Denote the comultiplication, the counit, and the
antipode by $\Delta$, $\ve$, and  by $S$, respectively.
We use the notions ``right comodule" and ``corepresentation" of $\AA$ as 
synonyms. By fixing a basis in the underlying vector space we identify
corepresentations and the corresponding matrices.  
Let $v$ (resp.\ $f$) be a corepresentation (resp.\  a representation) of
$\AA$. As usual $v^\fettc$ (resp.\  $f^\fettc$)  
denotes the contragredient corepresentation (resp.\  contragredient
representation) of $v$ (resp.\  of $f$). The
space of intertwiners of corepresentations $v$ and $w$ is  $\Mor(v,w)$.
We write $\Mor(v)$
for $\Mor(v,v)$. By $\End(V)$ and $V\ot W$ we always mean
$\End_\C(V)$ and $V\ot_\C W$, respectively. If $A$ is a linear mapping, $A^\fettt$
denotes the transpose of $A$ and $\tr A$ the trace of $A$.
Lower indices of $A$ always refer
to the components of a tensor product where $A$ acts (`leg numbering').
The unit matrix is denoted by $I$. 
Unless it is explicitly stated otherwise,
we use Einstein convention to sum over repeated indices.
 Set $\plus{a}=a-\ve(a)1$ for
$a\in\AA$ and $\plus{\AA}=\{\plus{a}|a\in\AA\}$.
We use Sweedler's notation for the coproduct $\Delta(a)=\sum a_{(1)}
\ot a_{(2)}$ and  for right
comodules  $\dr(\rho)=\sum \rho_{(0)}\ot \rho_{(1)}$.
The mapping $\Adr\colon\AA\to\AA\ot\AA$ defined by
 $\Adr a=\sum a_{(2)}\ot Sa_{(1)} a_{(3)}$ is called
 the {\em right adjoint coaction}\/ of $\AA$ on itself.
The mapping $b\rac a:=Sa_{(1)}ba_{(2)}$, $a\in \AA$, $b\in\BB$, is called the
{\em right adjoint action} of $\AA$ on $\BB$, where $\BB$ is an $\AA$-bimodule.

\section{Preliminaries}\label{prelim}
 In the next three  subsections we shall use the general framework of bicovariant differential calculus
developed by Woronowicz \cite{a-Woro2}, see also \cite[Chapter~14]{b-KS}. 
We collect the main notions and facts needed in what follows.

{\em Bicovariant bimodules.}
A {\em bicovariant bimodule}\/ over $\AA$ 
is a bimodule $\Gamm$ together with linear mappings 
$\dl\colon\Gamm\to\AA\ot\Gamm$ and $\dr\colon\Gamm\to\Gamm\ot\AA$ such that
$(\Gamm,\dl,\dr)$ is a bicomodule and 
$\dl(a\omega b)=\Delta(a)\dl(\omega)\Delta(b)$, 
and $\dr(a\omega b)=\Delta(a)\dr(\omega)\Delta(b)$
for $a,\,b\in\AA$ and $\omega\in\Gamm$.
Let $\Gamm$ be a bicovariant bimodule over $\AA$. 
We call the elements of the vector space 
$\Gamm_\ll=\{\om|\dl(\om)=1\ot\om\}$
(resp. $\Gamm_\rr=\{\om|\dr(\om)=\om\ot1\}$) {\em left-invariant} (resp. {\em right-invariant}\/). The elements of $\gi=\gl\cap\gr$ are called {\em bi-invariant}.
The structure of bicovariant bimodules has been  completely characterized 
by Theorems 2.3 and 2.4 in \cite{a-Woro2}. We recall the corresponding result:
Let $(\Gamm,\dl,\dr)$ be a bicovariant bimodule over $\AA$ and let 
$\{\om_i|i\in K\}$ be a finite linear basis of $\gl$. 
Then there exist matrices $v=(v^i_j)$
and $f=(f^i_j)$ of elements $v^i_j\in \AA$ and of functionals $f^i_j$ on
$\AA$, $i,j\in K$ such that $v$ is a matrix corepresentation, $f$ is matrix
representation of $\AA$, and
\begin{align}
\om_i\rac a&= f^i_n(a)\om_n,\label{e-rac}
\\
\dr(\om_i)&=\om_n\ot v^n_i,\label{e-rcac}
\end{align}
for  $a\in\AA$, $i\in K$. 
Conversely, if the  corepresentation $v$ and the representation $f$ satisfy
certain compatibility condition, then there exists a unique bicovariant
bimodule $\Gamm$ with \rf[e-rac] and \rf[e-rcac] and $\{\om_i|i\in K\}$ is a
basis of $\gl$. In this situation we simply  write $\Gamm=(v,f)$.

{\em Bicovariant first order differential calculi.}
A {\em first order differential calculus}\/ over $\AA$ (FODC for short)  is an
$\AA$-bimodule $\Gamm$ with  a linear mapping  $\dd\colon\AA\to\Gamm$ that
satisfies the Leibniz rule $\dd(ab)=\dd a{\cdot} b+a{\cdot}\dd b$ for
$a,\,b\in\AA$,  and $\Gamm$ is the linear span of elements $a\dd b$ with
$a,b\in\AA$.
\\
A FODC $\Gamm$  is called {\em bicovariant}\/ if there exist linear mappings
$\dl\colon\Gamm\to\AA\ot\Gamm$ and
$\dr\colon\Gamm\to\Gamm\ot\AA$ such that
\begin{align*}
\dl(a\dd b)&=\Delta(a)(\id\ot\dd)\Delta(b),
\\
\dr(a\dd b)&=\Delta(a)(\dd\ot\id)\Delta(b)
\end{align*}
for all $a,b\in\AA$. It turns out that $(\Gamm,\dl,\dr)$ is a bicovariant
bimodule. A bicovariant FODC is called {\em inner}\/ if there exists a bi-invariant
$1$-form   $\theta\in\Gamm$ such that
$$
\dd a=\theta a-a\theta,\quad a\in\AA. 
$$
By the dimension of a bicovariant FODC we mean the dimension of the vector
space $\gl$ of left-invariant 1-forms. 
Let $\Gamm$ be a bicovariant FODC over $\AA$. Then the set 
$$
\RR_\Gamm=\{a\in\plus{\AA}|\om(a)=0\}
$$
is an $\Adr$-invariant right ideal of $\plus{\AA}$. Here $\om\colon\AA\to\gl$
is the mapping 
\begin{align}\label{e-om}
\om(a)=Sa _{(1)}\dd a_{(2)}.
\end{align}
 Conversely, for any
$\Adr$-invariant right ideal $\RR$ of $\plus{\AA}$ there exists a bicovariant FODC
$\Gamm$ such that $\RR_\Gamm=\RR$, \cite[Proposition~14.7]{b-KS}.
\\
The linear space
$$
\XX_\Gamm=\{X\in\AA^\circ| X(1)=0 \text{ and } X(p)=0 \text{ for all } p\in
\RR_\Gamm\}
$$
is called the {\em quantum Lie algebra} of $\Gamm$.
We recall the main property. The space 
$\XX_\Gamm$ is an $\adr$-invariant subspace of the dual Hopf algebra
  $\AA^\circ$ satisfying $\Delta(X)-1\ot X\in\XX_\Gamm\ot \AA^\circ$ for
  $X\in\XX_\Gamm$, \cite[Corollary~14.10]{b-KS}.

{\em Higher order differential calculi.}
In this subsection we briefly repeat two concepts to construct higher order
differential calculi (DC for short)  for a given bicovariant FODC $\Gamm$. Let
$\Gamm=(v,f)$ be a bicovariant bimodule.
\\
Obviously the tensor product $\gten[k]=\Gamm\ota\cdots\ota\Gamm$ ($k$ factors)
is again a bicovariant bimodule. Define the tensor algebra
$\gten=\bigoplus_{k\ge0}\gten[k]$, $\gten[0]=\AA$, over $\AA$. This is also a bicovariant
bimodule. Since bicovariant bimodules are free left $\AA$-modules we always 
identify $(\Gamm\ota\cdots\ota\Gamm)_\ll$ and $\Gamm_\ll\ot\cdots\ot\Gamm_\ll$. 
This justifies our notation $\om_i\ot\om_j$ instead of $\om_i\ota\om_j$ for $\om_i,\om_j\in\gl$.
There exists  a unique isomorphism 
$\sig\colon\Gamm\ota\Gamm\to\Gamm\ota\Gamm$ of bicovariant bimodules
called the {\em braiding}\/   with
$\sig(\om\ott\rho)=\rho_{(0)}\ott(\om\rac\rho_{(1)})$,
$\om,\rho\in\Gamm_{\ll}$. Moreover $\sig$ fulfils the
braid equation 
$(\sig\ot\id)(\id\ot \sig)(\sig\ot \id)=(\id\ot \sig)(\sig\ot\id)(\id\ot\sig)$
in $\Gamm\ota\Gamm\ota\Gamm$.
Let  $\Js$ denote the two-sided ideal in $\gten$ generated by the kernel of
$A_2\colon\Gamm\ota\Gamm\to\Gamm\ota\Gamm$,  $A_2=\id-\sig$. We call
$\gds=\gten{/}\Js$ the {\em second antisymmetrizer exterior algebra} over
$\Gamm$. Since $\sig$ is a morphism of bicomodules, $(\Gamm\ota\Gamm)_\ll$ is
invariant under $\sig$. So there exist complex numbers $\sig^{ij}_{st}$ such
that $\sig(\om_s\ott\om_t)=\sig^{ij}_{st}\om_i\ott\om_j$. By
\cite[(3.15)]{a-Woro2} we have
\begin{equation}\label{fv}
\sig^{ij}_{st}=f^s_j(v^i_t).
\end{equation}
\\
Let $\SS\colon
\AA\to\gl\ot\gl$ be defined by
$$
\SS(a)=\om(a_{(1)})\ot\om(a_{(2)}).
$$
Let  $\Ju$ denote the two-sided ideal of $\gten$ generated by the vector
space $\SS(\RR_\Gamm)$. Then $\gdu=\gten{/}\Ju$ is called the {\em universal
  exterior algebra} over $\Gamm$. 
Both $\gds$ and $\gdu$ are $\NO$-graded algebras,  bicovariant bimodules over
$\AA$ as well as differential graded Hopf algebras over $\AA$. They are related by
$\Ju\subseteq\Js$. Their left-invariant subalgebras $\gdul$ and $\gdsl$ are
both quadratic algebras over the same vector space $\gl$. 

\section{Orthogonal and symplectic quantum groups, their standard FODC, and the main result}\label{osp}
In this section we recall general facts about orthogonal and symplectic quantum groups.
Throughout  $\AA$ denotes  
one of the Hopf algebras  $\ooq$ and $\ospq$ as
defined in \cite[Subsection\,1.4]{a-FRT}.
We give a brief introduction into the graphical calculus with morphisms of corepresentations of $\AA$ and we recall the construction of the standard bicovariant FODC over $\AA$. At the end we state our main result.
\\
  As usual we set  
$\eps=1$  in the orthogonal  and $\eps=-1$ in  the symplectic case.
Throughout  the deformation parameter $q$ is {\em not a root of unity},
 and $N\ge3$.
We always use the abbreviations  $\QM=q-q^{-1}$, $\qtwo=q+\qm$,
$\rt=\eps q^{N-\eps}$, and $\xx=1+\frac{\rt-\rtm}{q-q^{-1}}$. 
Recall that $R$ denotes the complex invertible $N^2\times N^2$-matrix
\cite[(1.9)]{a-FRT},
$\Rda[ab]{st}=R^{ba}_{st}$, and $C=(C^i_j)$,
$C^i_j=\eps_iq^{\vr_j}\delta_{ij'}$ defines the metric, see \cite[(1.10)]{a-FRT} for
details. The matrix $K$ is given by $K^{ab}_{st}=C^a_bB^s_t$, where $B=C^{-1}=\eps
C$. We need the diagonal matrix $D=B^\fettt C$.
Sometimes we use the notation $C^{ab}=C^a_b$, $C_{ab}=C^a_b$. Then
$(C^{ab})\in\End(\C,\C^N\ott\C^N)$ and $(C_{ab})\in\End(\C^N\ott\C^N,\C)$.
The  $N^2$ generators of $\AA$ are denoted
by  $u^i_j$, $i,\,j=1,\dots,N$, and we call $u=(u^i_j)$ the
fundamental matrix corepresentation. 
The element  $U=\sum_{i,j}D^j_iu^i_j$
is called the quantum trace.  
 Note that
\begin{equation}\label{mor}
(C^{ab})\in\Mor(1,u\ott u),\quad(C_{ab})\in\Mor(u\ott u,1),\quad C^\fettt
\in\Mor(u^\fettc,u)=\Mor(u,u^\fettc).
\end{equation}
For $T=(T^{ab}_{st})\in\End(\C^N\ott\C^N)$ define the $q$-trace $\trq
T\in\End(\C^N)$ by $(\trq T)^b_t=D^s_a T^{ab}_{st}$. We often use the following well known relations between $\Rda{}$, $\Rdam$, $K$,
and $D$. 
\begin{alignat}{2}
C^y_z\Rda[yz]{st}&=\rtm C^s_t,&\quad\Rda[ab]{yz}C^y_z&=\rtm C^a_b,\label{e-rc}
\\
\Rda{}-\Rdam&=\QM(I-K),&\label{e-rm}
\\
\xx&=\tr D,&\label{e-x}
\\
\trq\Rda{}&=\rt I,&\quad \trq I&=\xx I,\quad \trq K=I.\label{e-trq}
\end{alignat}
The mapping $g_i\mapsto \Rda{i,i+1}$, $e_i\mapsto K_{i,i+1}$ defines a
representation of the Birman-Wenzl-Murakami algebra $\mathrm{C}(q,\rt)$, \cite{a-Wenzl1}.
We shall give a brief introduction into the graphical calculus with morphisms,
see also \cite[Fig.\,1 and Fig.\,6]{a-SchSch1}. The calculus is justified in
\cite{a-Turaev1}. Using the graphical calculus formulas and proofs become more
transparent. In order to distinguish the places for the corepresentation $u$
and $u^\fettc$ we use  arrows in the graph. A vertex stands for $u$,
resp. $u^\fettc$, if the corresponding edge is downward directed, resp. upward
directed. Since for orthogonal and symplectic quantum groups $u$ and
$u^\fettc$ are isomorphic, it appears that one edge has {\em two}
directions. For instance, the intertwiner $C^\fettt\in\Mor(u,u^\fettc)$ is
represented by a vertical edge downward directed at the bottom and upward
directed at the top. Removing a curl by rotating part of the diagram clockwise (resp. anti-clockwise) acquires a factor $\rt$ (resp. $\rtm$) (First Reidemeister move). A closed loop gives the factor $\xx$. 
\begin{center}
\special{em:linewidth 0.4pt}
\unitlength 1.00mm
\linethickness{0.4pt}
\begin{picture}(145.33,53.67)
\put(17.67,36.67){\vector(1,-1){0.2}}
\emline{7.67}{46.67}{1}{17.67}{36.67}{2}
\put(7.67,37.00){\vector(-1,-1){0.2}}
\emline{12.00}{41.33}{3}{7.67}{37.00}{4}
\emline{18.00}{46.33}{5}{14.00}{42.33}{6}
\emline{8.00}{46.33}{7}{9.46}{48.08}{8}
\emline{9.46}{48.08}{9}{10.88}{49.25}{10}
\emline{10.88}{49.25}{11}{12.27}{49.84}{12}
\emline{12.27}{49.84}{13}{13.62}{49.85}{14}
\emline{13.62}{49.85}{15}{14.93}{49.29}{16}
\emline{14.93}{49.29}{17}{16.20}{48.14}{18}
\emline{16.20}{48.14}{19}{17.67}{46.00}{20}
\put(42.00,16.33){\vector(1,-1){0.2}}
\emline{32.00}{26.33}{21}{42.00}{16.33}{22}
\put(56.66,39.33){\vector(1,-1){0.2}}
\emline{46.66}{49.33}{23}{56.66}{39.33}{24}
\put(95.66,37.67){\vector(1,-1){0.2}}
\emline{85.66}{47.67}{25}{95.66}{37.67}{26}
\put(32.00,16.67){\vector(-1,-1){0.2}}
\emline{36.33}{21.00}{27}{32.00}{16.67}{28}
\put(46.66,39.67){\vector(-1,-1){0.2}}
\emline{51.00}{44.00}{29}{46.66}{39.67}{30}
\put(85.66,38.00){\vector(-1,-1){0.2}}
\emline{90.00}{42.33}{31}{85.66}{38.00}{32}
\emline{42.33}{26.00}{33}{38.33}{22.00}{34}
\emline{57.00}{49.00}{35}{53.00}{45.00}{36}
\emline{96.00}{47.33}{37}{92.00}{43.33}{38}
\emline{32.00}{26.33}{39}{30.17}{25.25}{40}
\emline{30.17}{25.25}{41}{28.90}{24.13}{42}
\emline{28.90}{24.13}{43}{28.18}{22.96}{44}
\emline{28.18}{22.96}{45}{28.01}{21.75}{46}
\emline{28.01}{21.75}{47}{28.40}{20.49}{48}
\emline{28.40}{20.49}{49}{29.34}{19.19}{50}
\emline{29.34}{19.19}{51}{32.00}{17.00}{52}
\put(21.00,42.00){\makebox(0,0)[cc]{$=\rtm$}}
\put(40.33,40.66){\makebox(0,0)[cc]{,}}
\emline{46.33}{39.67}{53}{47.54}{37.65}{54}
\emline{47.54}{37.65}{55}{48.75}{36.14}{56}
\emline{48.75}{36.14}{57}{49.97}{35.14}{58}
\emline{49.97}{35.14}{59}{51.18}{34.64}{60}
\emline{51.18}{34.64}{61}{52.39}{34.65}{62}
\emline{52.39}{34.65}{63}{53.60}{35.17}{64}
\emline{53.60}{35.17}{65}{54.82}{36.20}{66}
\emline{54.82}{36.20}{67}{57.00}{39.33}{68}
\put(61.33,42.33){\makebox(0,0)[cc]{$=\rtm$}}
\put(80.33,40.33){\makebox(0,0)[cc]{,}}
\put(100.99,37.67){\vector(-1,-1){0.2}}
\emline{111.33}{48.00}{69}{100.99}{37.67}{70}
\put(110.66,37.33){\vector(1,-1){0.2}}
\emline{106.33}{41.67}{71}{110.66}{37.33}{72}
\emline{100.66}{47.67}{73}{104.99}{43.33}{74}
\put(99.33,42.67){\makebox(0,0)[cc]{$-$}}
\put(113.33,42.67){\makebox(0,0)[cc]{$=\QM$}}
\put(116.99,37.33){\vector(0,-1){0.2}}
\emline{116.99}{47.67}{75}{116.99}{37.33}{76}
\put(122.33,37.33){\vector(0,-1){0.2}}
\emline{122.33}{47.67}{77}{122.33}{37.33}{78}
\put(126.32,43.00){\makebox(0,0)[cc]{$-\QM$}}
\put(45.00,21.67){\makebox(0,0)[cc]{$=\rt$}}
\put(49.33,16.67){\vector(0,-1){0.2}}
\emline{49.33}{26.33}{79}{49.33}{16.67}{80}
\put(53.00,20.33){\makebox(0,0)[cc]{,}}
\put(78.67,22.33){\makebox(0,0)[cc]{$=\xx$}}
\put(84.33,17.00){\vector(0,-1){0.2}}
\emline{84.33}{26.00}{81}{84.33}{17.00}{82}
\put(88.00,21.00){\makebox(0,0)[cc]{,}}
\put(118.33,21.67){\makebox(0,0)[cc]{$=$}}
\put(123.66,16.00){\vector(0,-1){0.2}}
\emline{123.66}{26.00}{83}{123.66}{16.00}{84}
\put(127.66,20.66){\makebox(0,0)[cc]{.}}
\put(74.33,17.00){\vector(0,-1){0.2}}
\emline{74.33}{26.00}{85}{74.33}{17.00}{86}
\put(66.66,17.00){\vector(0,-1){0.2}}
\emline{66.66}{26.00}{87}{66.66}{17.00}{88}
\emline{66.67}{26.33}{89}{64.84}{25.25}{90}
\emline{64.84}{25.25}{91}{63.57}{24.13}{92}
\emline{63.57}{24.13}{93}{62.85}{22.96}{94}
\emline{62.85}{22.96}{95}{62.68}{21.75}{96}
\emline{62.68}{21.75}{97}{63.07}{20.49}{98}
\emline{63.07}{20.49}{99}{64.01}{19.19}{100}
\emline{64.01}{19.19}{101}{66.67}{17.00}{102}
\emline{104.00}{26.33}{103}{102.18}{25.25}{104}
\emline{102.18}{25.25}{105}{100.90}{24.13}{106}
\emline{100.90}{24.13}{107}{100.18}{22.96}{108}
\emline{100.18}{22.96}{109}{100.01}{21.75}{110}
\emline{100.01}{21.75}{111}{100.40}{20.49}{112}
\emline{100.40}{20.49}{113}{101.34}{19.19}{114}
\emline{101.34}{19.19}{115}{104.00}{17.00}{116}
\put(24.66,36.67){\vector(-1,-1){0.2}}
\emline{26.66}{38.67}{117}{24.66}{36.67}{118}
\put(35.00,36.67){\vector(3,-4){0.2}}
\emline{26.33}{38.33}{119}{27.68}{39.65}{120}
\emline{27.68}{39.65}{121}{29.12}{40.26}{122}
\emline{29.12}{40.26}{123}{30.64}{40.18}{124}
\emline{30.64}{40.18}{125}{32.25}{39.39}{126}
\emline{32.25}{39.39}{127}{35.00}{36.67}{128}
\put(132.33,37.33){\vector(-1,-1){0.2}}
\emline{134.33}{39.33}{129}{132.33}{37.33}{130}
\put(142.66,37.33){\vector(3,-4){0.2}}
\emline{134.00}{39.00}{131}{135.35}{40.32}{132}
\emline{135.35}{40.32}{133}{136.79}{40.93}{134}
\emline{136.79}{40.93}{135}{138.31}{40.84}{136}
\emline{138.31}{40.84}{137}{139.91}{40.05}{138}
\emline{139.91}{40.05}{139}{142.66}{37.33}{140}
\put(104.00,16.67){\vector(-1,-1){0.2}}
\emline{106.00}{18.67}{141}{104.00}{16.67}{142}
\put(114.33,16.67){\vector(3,-4){0.2}}
\emline{105.66}{18.33}{143}{107.02}{19.65}{144}
\emline{107.02}{19.65}{145}{108.45}{20.26}{146}
\emline{108.45}{20.26}{147}{109.98}{20.18}{148}
\emline{109.98}{20.18}{149}{111.58}{19.39}{150}
\emline{111.58}{19.39}{151}{114.33}{16.67}{152}
\put(67.67,45.33){\vector(1,-1){0.2}}
\emline{66.00}{47.00}{153}{67.67}{45.33}{154}
\emline{67.67}{45.33}{155}{68.94}{43.81}{156}
\emline{68.94}{43.81}{157}{70.42}{43.25}{158}
\emline{70.42}{43.25}{159}{72.10}{43.64}{160}
\emline{72.10}{43.64}{161}{74.00}{45.00}{162}
\put(74.00,44.67){\vector(-4,-3){0.2}}
\emline{76.34}{46.67}{163}{74.00}{44.67}{164}
\put(134.33,46.33){\vector(1,-1){0.2}}
\emline{132.66}{48.00}{165}{134.33}{46.33}{166}
\emline{134.33}{46.33}{167}{135.60}{44.81}{168}
\emline{135.60}{44.81}{169}{137.08}{44.25}{170}
\emline{137.08}{44.25}{171}{138.76}{44.64}{172}
\emline{138.76}{44.64}{173}{140.66}{46.00}{174}
\put(140.66,45.67){\vector(-4,-3){0.2}}
\emline{143.00}{47.67}{175}{140.66}{45.67}{176}
\put(105.66,24.67){\vector(1,-1){0.2}}
\emline{103.99}{26.34}{177}{105.66}{24.67}{178}
\emline{105.66}{24.67}{179}{106.93}{23.15}{180}
\emline{106.93}{23.15}{181}{108.41}{22.58}{182}
\emline{108.41}{22.58}{183}{110.10}{22.98}{184}
\emline{110.10}{22.98}{185}{111.99}{24.34}{186}
\put(111.99,24.01){\vector(-4,-3){0.2}}
\emline{114.33}{26.01}{187}{111.99}{24.01}{188}
\put(145.33,41.67){\makebox(0,0)[cc]{,}}
\end{picture}
\\
Figure\,1: The graphical representation of \rf[e-rc], \rf[e-rm], and \rf[e-trq].
\end{center}
The matrix $\Rda{}$ has the spectral decomposition
$$
\Rda{}=q\pda[+]-\qm \pda[-]+\rtm \pda[0],
$$
where $\pda$, $\tau\in\{+,-,0\}$, is idempotent.
\\
We recall  the method of Jur\v co \cite{a-Jurco1} to construct bicovariant FODC
over $\AA$. For the more general construction of bicovariant FODC over
coquasitriangular Hopf algebras see \cite[Section~14.5]{b-KS}.
Let  $\ell^\pm=(\ell^{\pm i}_j)$ be the $N\times N$-matrix of linear
functionals $\ell^{\pm i}_j$ on $\AA$ as defined in
in \cite[Section\,2]{a-FRT}. Recall that $\ell^\pm$ is uniquely determined by
${\ell^{\pm}}^i_j(u^m_n)=(\Rda[\pm1]{})_{nj}^{im}$ and the
property that  $\ell^\pm\colon\AA\to \End(\C^N)$ is a unital algebra
homomorphism. Note that $\ell^{\pm i}_j(S
u^m_n)=(\Rda[\mp 1]{})^{mi}_{jn}$. Define the bicovariant bimodules 
\begin{align*}
\Gamm_\pm=(u\ot u^\fettc,\vepm\ot\ell^{-\fettc}\ot\ell^{+}),
\end{align*}
where $\ve_+=\ve$ and  $\vem$ is the character on $\AA$ given by $\vem(u^i_j)=-\delta_{ij}$.
The structure of $\gpm$ can easily be described as follows. There exists a
basis $\{\theta_{ij}|i,\,j=1,\dots,N\}$ of $(\gpm)_\ll$ such that  the right
adjoint 
action and the right coaction  are given by
\begin{align*}
\theta_{ij}\rac a&=\vepm(a_{(1)})S(\ell^{-m}_i)\ell^{+j}_n(a_{(2)})\theta_{mn},\quad a\in \AA,
\\         
\dr \theta_{ij}&=\theta_{mn}\ot u^m_i(u^\fettc)^n_j,\quad i,\,j=1,\dots,N.
\end{align*}
In particular
\begin{equation}\label{e-thetau}
\begin{split}
\theta_{ij}\rac u^s_t&=\pm\Rda[sm]{iy}\Rda[jy]{tn}\theta_{mn},
\\
\theta\rac u^s_t     &=\pm(\Rda[2]{})^{sm}_{tn}\theta_{mn},
\end{split}
\end{equation}
where $\theta=\sum_i\theta_{ii}$ is the unique up to scalars
bi-invariant element. Defining 
\begin{equation}
\label{inner}
\dd a=\theta a-a\theta
\end{equation}
for $a\in \AA$, $(\gpm,\dd)$ becomes a bicovariant FODC over $\AA$. 
The  basis $\{X^\pm_{ij}\}$ of the quantum Lie algebra
$\XX_\pm$ dual to $\{\theta_{ij}\}$ is given by
\begin{align*}
X^\pm_{ij}:=\ve_\pm\ell^i_j-\delta_{ij}:=\ve_\pm S(\ell^{-i}_y)\ell^{+y}_j-\delta_{ij}.
\end{align*}
One easily checks that $X^\pm_0:=D^j_iX^\pm_{ij}=\ve_\pm D^j_i\ell^i_j-\xx$
is an $\adr$-invariant element of $\AA^\circ$. The braiding $\sig$ of $\gpm$ can be obtained as follows. Inserting
$v=u\ott u^\fettc$ and \mbox{$f_\pm=\vepm\ott\ell^{-\fettc}\ott\ell^{+}$} into
equation \rf[fv] the braiding matrices of $\Gamm_+$ and $\Gamm_-$ coincide
\begin{equation}\label{sig}
\sig=\Rli[-]{23}\Rda{12}\Rchm[34]\Rre{23},
\end{equation}
where the  matrices $\Rre{}$, $\Rch{}$, and $\Rli[-]{}$
 are
defined as follows. For a complex $N^2{\times}N^2$-matrix $T$ with 
$\hat{T}\in\Mor(u\ott u)$ define the matrices $\check{T}^{ab}_{st}=\hat{T}^
{ts}_{ba}$, $\acute{T}^{ab}_{st}=\hat{T}^{sa}_{tb}$,
and $\grave{T^-}=(\acute{T})^{-1}$. Note that $\acute{T}\in\Mor(u^\fettc\ott 
u,u\ott u^\fettc)$ and $\check{T}\in\Mor(u^\fettc\ott u^\fettc)$.
\\
Now we can formulate our main result.
\begin{thm}\label{t-1}
Let $\AA$ be one of the Hopf algebras $\ooq$ or $\ospq$,  \mbox{$N\ge3$}, and  $q$
not a root of unity. Let  $\Gamm$ be one of the
bicovariant FODC $\Gamm_\pm$, and
$2\xx+{(q-\qm)}(\rt-\rtm)\ne0$ in case of $\Gamm_-$.
Denote the  unique up to scalars bi-invariant $1$-form by $\theta$.
\\
{\em (i)} Then the quotient  $\gdu{/}(\theta^2)$  and
the second antisymmetrizer algebra $\gds$ are isomorphic bicovariant bimodules.
\\
{\em (ii)} The bi-invariant  $2$-form $\theta^2$ is central in $\gdu$. The
calculus $\gdu$ is  inner, 
i.\,e. 
\begin{align}\label{e-inner}
\dd\rho=\theta\land\rho-(-1)^n\rho\land\theta,\quad \rho\in\gdu[n].
\end{align}
\end{thm}
\begin{remark}
Theorem~1 is true for the quantum group $\slqn[2]$ and the $4D_\pm$ bicovariant
FODC  as well, \cite[Theorem~3.3\,(iii)]{a-Sch1}. In cases $\slqn$
and $\glqn$, $N\ge3$, we have $\gdu\cong\gds$, \cite[Theorem~3.3\,(ii)]{a-Sch1}. 
For the quantum super group $\mathrm{GL}_q(m|n)$ the relation $\gdu\cong\gds$ was proved in \cite[Section\,5.3]{a-LyuSud1}.
\end{remark}
\begin{remark}
The isomorphism of bicovariant bimodules $\gdu{/}(\theta^2)$ and $\gds$
implies its isomorphy as differential graded Hopf algebras.
\end{remark}
\section{Proof of the Theorem}\label{proof}
In the first part of this section we study the duality of
$\Gamm_\ll\ott\Gamm_\ll$ and $\XX\ott\XX$ in more detail. In the second part 
we examine how $\Gamm\ota\Gamm$ splits into bicovariant subbimodules. We shall
prove that the space of bi-invariant elements of  $\Gamm\ota\Gamm$ 
 generates the whole bimodule $\ker A_2$.

{\em Duality.}
There is a useful criterion  to describe the dimension of the 
space of left-invariant \mbox{2-forms} of $\gdu$ in terms of the quantum Lie algebra. 
\begin{lemma}\label{l-1}
Let  $\AA$ be an arbitrary Hopf algebra, $\Gamm$  a left-covariant FODC over
$\AA$  with quantum Lie algebra $\XX$, and  $\gdu$ the
universal differential calculus over $\Gamm$. Then 
\begin{align*}
\dim\gdul[2]=\dim\{T\in\XX\ot\XX|\, \mu(T)\in\XX\},
\end{align*}
where $\mu\colon \XX\ot\XX\to\AA^\circ$ denotes the multiplication map.
\end{lemma}
\begin{proof}  We use the following simple lemma from linear algebra without
proof.
Let $B\colon V\times W\to \C$ be a non-degenerate linear pairing of
finite dimensional 
vector spaces and $U$ a subspace of $V$. Then the induced pairing
$\ov{B}\colon V/U\times U^\perp\to\C$ with $U^\perp=\{w\in W\,|\, B(u,w)=0
\text{ for } u\in U\}$ is also non-degenerate. Applying this lemma to the
non-degenerate pairing $\pair\colon\gl\ott\gl\times\XX\ott\XX\to \C$,
 \cite[p.~164]{a-Woro2}, and $U=\SS(\RR)$  we have  
$T=\alpha^{ij}X_i\ott X_j\in U^\perp$ if and only if
\begin{align*}
0=\langle\om(p_{(1)})\ott\om(p_{(2)}),\,\,\alpha^{ij}X_i\ott X_j\rangle=\alpha^{ij}X_i(p_{(1)})X_j(p_{(2)})=\mu(T)(p)
\end{align*}
for $p\in \RR$. Hence $T\in U^\perp$ if and only if
$\mu(T)\in\XX$. Consequently  $U^\perp=\mu^{-1}(\XX)$, where $\mu^{-1}(\XX)$
denotes the pre-image of $\XX$ under $\mu$.
Since the induced pairing is also non-degenerate  and
$\gdul[2]=\gl\ott\gl{/}\SS(\RR)$ by definition, the assertion of the lemma is proved.
\end{proof}
\begin{remark}
Suppose $\Gamm$ to be bicovariant. Since for 
$f\colon V\to V$  linear, $(\ker f)^\perp=\im f^\fettt$,  the  pairing also factorizes to a non-degenerate pairing of $\gdsl[2]\times
\XX{\land}\XX$, where $\XX{\land}\XX=A_2^\fettt(\XX\ott\XX)$ and $A_2^\fettt$ is
the dual mapping to $A_2\uhr(\gl\ot\gl)$. 
\end{remark}
We proceed with a result for a dual pairing of  a comodule and a module.
\begin{prp}\label{p-1}
Let $V$ be  a right $\AA$-comodule, $W$ a right $\AA^\circ$-module, and 
$\langle\,,\,\rangle\colon V{\times}W\to\C$ a non-degenerate dual pairing of
vector spaces. Suppose further that 
\begin{align*}
\langle v,w\cdot f\rangle=\langle v_{(0)}f(v_{(1)}), w\rangle,
\end{align*}
for $v\in V$, $w\in W$, and $f\in\AA^\circ$. 
\\
If $P\in\Mor(V)$ then $P^\fettt\in\Mor(W)$. If in addition $P^2=P$, then  the
induced pairing $\im P\times \im P^\fettt \to \C$ is also non-degenerate.
\end{prp}
\begin{proof}
Since $P\in\Mor(V)$, $Pv_{(0)}\ot v_{(1)}=(Pv)_{(0)}\ot (Pv)_{(1)}$. For $v\in
V$, $w\in W$, and $f\in \AA^\circ$ we thus get
\begin{align*}
\langle v,P^\fettt w\cdot f\rangle&=\langle v_{(0)}f(v_{(1)}),P^\fettt
w\rangle
\\
&=\langle P v_{(0)}f(v_{(1)}),w\rangle
\\
&=\langle (Pv)_{(0)}f((Pv)_{(1)}),w\rangle
\\
&=\langle Pv,w\cdot f\rangle
\\
&=\langle v, P^\fettt (w\cdot f)\rangle.
\end{align*}
Since the pairing is non-degenerate the first assertion follows.
\\
Since $P$ and $P^\fettt$ are morphisms, the corresponding subspaces are
invariant. Let $v_0\in \im P$, i.\,e. $v_0=Pv_0$,  and suppose $0=\langle v_0,P^\fettt w\rangle$ for
all $w\in W$. Then $0=\langle Pv_0, P^\fettt w\rangle=\langle
P^2v_0,w\rangle=\langle v_0,w\rangle$. Since the pairing is non-degenerate,
$v_0=0$. Similarly one shows that $\im P$ separates the elements of $\im P^\fettt$.
\end{proof}
\begin{cor} \label{c-1} Let  $\AA$ be an arbitrary Hopf algebra, $\Gamm$  a
  bicovariant FODC over $\AA$  with quantum Lie algebra $\XX$, and $P\in\Mor(\dr)$, $P^2=P$. We restrict $\dr$ to $(\Gamm\ota\Gamm)_\ll$ or a suitable quotient.
\\
Then  $\im P$ is a  $\dr$-invariant subspace of  $(\Gamm\ota\Gamm)_\ll$ ($\gdul[2]$
 resp.{} $\gdsl[2]$), and  $\im P^\fettt$ is an
$\adr$-invariant subspace of  $\XX\ott\XX$ ($\mu^{-1}(\XX)$ resp.{}
$\XX{\land}\XX$). The induced pairing $\im P{\times}\im P^\fettt\to\C$ is non-degenerate.
\end{cor}
\begin{proof}
(i) Since $\SS(\RR)$ and  $(\ker A_2)_\ll$ are $\dr$-invariant, and since 
$\mu^{-1}(\XX)$ and $A_2^\fettt(\XX\ott\XX)$ are $\adr$-invariant, the  mappings $\dr$ 
and $\adr$ are  well-defined on both quotients $\gdul[2]$ and $\gdsl[2]$
resp. $\mu^{-1}(\XX\ott\XX)$ and $\XX{\land}\XX$.
\\
 It follows from 
\cite[(5.17) and (5.21)]{a-Woro2} that for $\rho\in(\Gamm\ota\Gamm)_\ll$,
$Y\in\XX\ott\XX$, and $f\in\AA^\circ$
$$
\langle \rho_{(0)}f(\rho_{(1)})\,,\,\,Y\rangle=\langle\rho\,,\,Y\rac f\rangle.
$$
Thus Proposition\,\ref{p-1} applies to our situation.
\end{proof}
Our next aim is to compare the bi-invariant components of $\Gamm\ota\Gamm$,
$\gds[2]$, and $\gdu[2]$ resp. of $\XX\ott\XX$, $\XX{\land}\XX$, and $\mu^{-1}(\XX)$. Let $\BB$ be a right
$\AA^\circ$-module with respect to $\adr$. For the space of invariants we use
the notation $\BB_0=\{b\in\BB| b\rac f=\ve(f)b, \, f\in\AA^\circ\}$.
\begin{lemma}\label{l-biinv}
Let $\AA$ be one of the Hopf algebras $\ooq$ or $\ospq$,  $N\ge3$, $\Gamm$ one of the $N^2$-dimensional bicovariant FODC $\Gamm_\pm$ over
$\AA$ and let $\XX$ be the corresponding quantum Lie algebra. Then we have
\begin{alignat*}{4}
\text{{\rm (i)}}&&\qquad \dim \gteni[2]&=3,&\quad \dim(\XX\ott\XX)_0&=3,
\\
\text{{\rm (ii)}}&&\qquad \dim \gdsi[2]&=0,&\quad \dim(\XX{\land}\XX)_0&=0,
&\quad \dim(\ker A_2)_\ii&=3,
\\
\text{{\rm (iii)}}&&\qquad \dim \gdui[2]&=1,&\quad \mu^{-1}(\XX)_0&=
\langle T\rangle,
&\quad \dim(\SS(\RR)_\ii)&=2, 
\end{alignat*}
where 
$$
T=X_{ij}\ot X_{mn}\,B^i_y\Rda[jy]{mz}C^n_z.
$$
\end{lemma}
\begin{proof}
(i)  It is well known that $\dim\Mor(u\ott u)=3$, and $I,\,\Rda{}$, and $K$
form a linear basis of $\Mor(u\ott u)$. Using \rf[mor] it is easy to see that the
mapping $E\mapsto (B^s_zC^r_yE^{ab}_{yz})$ defines a linear  isomorphism
$\Mor(u\ott u)\to \Mor(1,u\ott u^\fettc\ott u\ott u^\fettc)$. Since
$\rho=\alpha^{ijmn}\theta_{ij}\ott\theta_{mn}\in\Gamm\ota\Gamm$ is
bi-invariant if and only if $(\alpha^{ijmn})\in\Mor(1,u\ott u^\fettc\ott u\ott
u^\fettc)$,  (i)  is proved. 
\\
(ii) The elements $\theta\ott\theta$, 
$\eta=D^k_j\theta_{ik}\ott\theta_{ji}$, and $\xi=C^i_z\Rda[ym]{zn}B^y_j\theta_{ij}\ott\theta_{mn}$ form a basis of $\gteni[2]$.
Using the graphical calculus it is not difficult to
check that $\sig$ acts as the identity on $\gteni[2]$:
\begin{center}
\special{em:linewidth 0.4pt}
\unitlength 1.00mm
\linethickness{0.4pt}
\begin{picture}(152.67,25.67)
\put(15.00,25.33){\vector(-1,1){0.2}}
\emline{19.00}{21.33}{1}{15.00}{25.33}{2}
\emline{10.33}{25.67}{3}{17.00}{19.00}{4}
\emline{20.33}{25.00}{5}{18.67}{23.00}{6}
\emline{18.67}{23.00}{7}{18.67}{23.00}{8}
\emline{18.67}{23.00}{9}{18.67}{23.00}{10}
\emline{21.00}{19.33}{11}{27.67}{12.67}{12}
\put(23.33,12.67){\vector(1,-1){0.2}}
\emline{19.00}{16.67}{13}{23.33}{12.67}{14}
\put(26.00,25.67){\vector(1,1){0.2}}
\emline{13.33}{13.00}{15}{26.00}{25.67}{16}
\put(8.33,12.67){\vector(-1,-1){0.2}}
\emline{14.67}{19.00}{17}{8.33}{12.67}{18}
\put(58.66,25.00){\vector(-1,1){0.2}}
\emline{62.66}{21.00}{19}{58.66}{25.00}{20}
\emline{54.00}{25.34}{21}{60.66}{18.67}{22}
\emline{64.00}{24.67}{23}{62.33}{22.67}{24}
\emline{62.33}{22.67}{25}{62.33}{22.67}{26}
\emline{62.33}{22.67}{27}{62.33}{22.67}{28}
\emline{64.66}{19.00}{29}{71.33}{12.34}{30}
\put(67.00,12.34){\vector(1,-1){0.2}}
\emline{62.66}{16.34}{31}{67.00}{12.34}{32}
\put(69.66,25.34){\vector(1,1){0.2}}
\emline{57.00}{12.67}{33}{69.66}{25.34}{34}
\put(52.00,12.34){\vector(-1,-1){0.2}}
\emline{58.33}{18.67}{35}{52.00}{12.34}{36}
\put(108.67,25.00){\vector(-1,1){0.2}}
\emline{112.67}{21.00}{37}{108.67}{25.00}{38}
\emline{104.00}{25.33}{39}{110.67}{18.66}{40}
\emline{114.00}{24.66}{41}{112.33}{22.66}{42}
\emline{112.33}{22.66}{43}{112.33}{22.66}{44}
\emline{112.33}{22.66}{45}{112.33}{22.66}{46}
\emline{114.67}{19.00}{47}{121.33}{12.33}{48}
\put(117.00,12.33){\vector(1,-1){0.2}}
\emline{112.67}{16.33}{49}{117.00}{12.33}{50}
\put(119.67,25.33){\vector(1,1){0.2}}
\emline{107.00}{12.66}{51}{119.67}{25.33}{52}
\put(102.00,12.33){\vector(-1,-1){0.2}}
\emline{108.33}{18.66}{53}{102.00}{12.33}{54}
\put(27.67,12.67){\vector(1,4){0.2}}
\emline{23.33}{12.67}{55}{23.89}{10.60}{56}
\emline{23.89}{10.60}{57}{24.47}{9.19}{58}
\emline{24.47}{9.19}{59}{25.06}{8.46}{60}
\emline{25.06}{8.46}{61}{25.66}{8.39}{62}
\emline{25.66}{8.39}{63}{26.27}{8.99}{64}
\emline{26.27}{8.99}{65}{26.90}{10.26}{66}
\emline{26.90}{10.26}{67}{27.67}{12.67}{68}
\put(13.33,12.33){\vector(1,2){0.2}}
\emline{8.33}{12.67}{69}{8.99}{10.59}{70}
\emline{8.99}{10.59}{71}{9.65}{9.17}{72}
\emline{9.65}{9.17}{73}{10.34}{8.40}{74}
\emline{10.34}{8.40}{75}{11.03}{8.28}{76}
\emline{11.03}{8.28}{77}{11.73}{8.83}{78}
\emline{11.73}{8.83}{79}{12.45}{10.03}{80}
\emline{12.45}{10.03}{81}{13.33}{12.33}{82}
\put(71.33,12.00){\vector(1,2){0.2}}
\emline{52.00}{12.00}{83}{53.51}{10.02}{84}
\emline{53.51}{10.02}{85}{55.00}{8.33}{86}
\emline{55.00}{8.33}{87}{56.47}{6.94}{88}
\emline{56.47}{6.94}{89}{57.91}{5.85}{90}
\emline{57.91}{5.85}{91}{59.32}{5.05}{92}
\emline{59.32}{5.05}{93}{60.70}{4.55}{94}
\emline{60.70}{4.55}{95}{62.06}{4.34}{96}
\emline{62.06}{4.34}{97}{63.40}{4.43}{98}
\emline{63.40}{4.43}{99}{64.71}{4.81}{100}
\emline{64.71}{4.81}{101}{65.99}{5.49}{102}
\emline{65.99}{5.49}{103}{67.25}{6.47}{104}
\emline{67.25}{6.47}{105}{68.48}{7.74}{106}
\emline{68.48}{7.74}{107}{69.69}{9.31}{108}
\emline{69.69}{9.31}{109}{71.33}{12.00}{110}
\put(56.66,12.34){\vector(-3,4){0.2}}
\emline{67.00}{12.34}{111}{65.96}{10.44}{112}
\emline{65.96}{10.44}{113}{64.85}{9.09}{114}
\emline{64.85}{9.09}{115}{63.67}{8.28}{116}
\emline{63.67}{8.28}{117}{61.09}{8.28}{118}
\emline{61.09}{8.28}{119}{59.68}{9.09}{120}
\emline{59.68}{9.09}{121}{58.21}{10.44}{122}
\emline{58.21}{10.44}{123}{56.66}{12.34}{124}
\put(107.00,12.00){\vector(-3,4){0.2}}
\emline{122.00}{11.66}{125}{120.58}{9.81}{126}
\emline{120.58}{9.81}{127}{119.17}{8.35}{128}
\emline{119.17}{8.35}{129}{117.77}{7.27}{130}
\emline{117.77}{7.27}{131}{116.37}{6.58}{132}
\emline{116.37}{6.58}{133}{114.97}{6.26}{134}
\emline{114.97}{6.26}{135}{113.58}{6.33}{136}
\emline{113.58}{6.33}{137}{112.20}{6.79}{138}
\emline{112.20}{6.79}{139}{110.82}{7.62}{140}
\emline{110.82}{7.62}{141}{109.45}{8.84}{142}
\emline{109.45}{8.84}{143}{107.00}{12.00}{144}
\emline{101.67}{11.66}{145}{102.12}{9.42}{146}
\emline{102.12}{9.42}{147}{102.99}{7.68}{148}
\emline{102.99}{7.68}{149}{104.29}{6.46}{150}
\emline{104.29}{6.46}{151}{106.02}{5.74}{152}
\emline{106.02}{5.74}{153}{109.67}{5.66}{154}
\emline{117.33}{11.66}{155}{115.47}{9.51}{156}
\emline{115.47}{9.51}{157}{112.33}{8.00}{158}
\put(47.67,18.67){\vector(1,4){0.2}}
\emline{43.33}{18.67}{159}{43.89}{16.60}{160}
\emline{43.89}{16.60}{161}{44.47}{15.19}{162}
\emline{44.47}{15.19}{163}{45.06}{14.46}{164}
\emline{45.06}{14.46}{165}{45.66}{14.39}{166}
\emline{45.66}{14.39}{167}{46.27}{14.99}{168}
\emline{46.27}{14.99}{169}{46.90}{16.26}{170}
\emline{46.90}{16.26}{171}{47.67}{18.67}{172}
\put(39.00,18.33){\vector(1,2){0.2}}
\emline{34.00}{18.67}{173}{34.66}{16.59}{174}
\emline{34.66}{16.59}{175}{35.32}{15.17}{176}
\emline{35.32}{15.17}{177}{36.01}{14.40}{178}
\emline{36.01}{14.40}{179}{36.70}{14.28}{180}
\emline{36.70}{14.28}{181}{37.40}{14.83}{182}
\emline{37.40}{14.83}{183}{38.12}{16.03}{184}
\emline{38.12}{16.03}{185}{39.00}{18.33}{186}
\put(31.00,18.00){\makebox(0,0)[cc]{$=$}}
\put(96.33,19.00){\vector(1,2){0.2}}
\emline{77.00}{19.00}{187}{78.51}{17.02}{188}
\emline{78.51}{17.02}{189}{80.00}{15.33}{190}
\emline{80.00}{15.33}{191}{81.47}{13.94}{192}
\emline{81.47}{13.94}{193}{82.91}{12.85}{194}
\emline{82.91}{12.85}{195}{84.32}{12.05}{196}
\emline{84.32}{12.05}{197}{85.70}{11.55}{198}
\emline{85.70}{11.55}{199}{87.06}{11.34}{200}
\emline{87.06}{11.34}{201}{88.40}{11.43}{202}
\emline{88.40}{11.43}{203}{89.71}{11.81}{204}
\emline{89.71}{11.81}{205}{90.99}{12.49}{206}
\emline{90.99}{12.49}{207}{92.25}{13.47}{208}
\emline{92.25}{13.47}{209}{93.48}{14.74}{210}
\emline{93.48}{14.74}{211}{94.69}{16.31}{212}
\emline{94.69}{16.31}{213}{96.33}{19.00}{214}
\put(81.66,19.34){\vector(-3,4){0.2}}
\emline{92.00}{19.34}{215}{90.96}{17.44}{216}
\emline{90.96}{17.44}{217}{89.85}{16.09}{218}
\emline{89.85}{16.09}{219}{88.67}{15.28}{220}
\emline{88.67}{15.28}{221}{86.09}{15.28}{222}
\emline{86.09}{15.28}{223}{84.68}{16.09}{224}
\emline{84.68}{16.09}{225}{83.21}{17.44}{226}
\emline{83.21}{17.44}{227}{81.66}{19.34}{228}
\put(74.33,19.00){\makebox(0,0)[cc]{$=$}}
\put(135.00,17.34){\vector(-3,4){0.2}}
\emline{150.00}{17.00}{229}{148.58}{15.15}{230}
\emline{148.58}{15.15}{231}{147.17}{13.69}{232}
\emline{147.17}{13.69}{233}{145.77}{12.61}{234}
\emline{145.77}{12.61}{235}{144.37}{11.92}{236}
\emline{144.37}{11.92}{237}{142.97}{11.60}{238}
\emline{142.97}{11.60}{239}{141.58}{11.67}{240}
\emline{141.58}{11.67}{241}{140.20}{12.13}{242}
\emline{140.20}{12.13}{243}{138.82}{12.96}{244}
\emline{138.82}{12.96}{245}{137.45}{14.18}{246}
\emline{137.45}{14.18}{247}{135.00}{17.34}{248}
\emline{129.67}{17.00}{249}{130.12}{14.76}{250}
\emline{130.12}{14.76}{251}{130.99}{13.02}{252}
\emline{130.99}{13.02}{253}{132.29}{11.80}{254}
\emline{132.29}{11.80}{255}{134.02}{11.08}{256}
\emline{134.02}{11.08}{257}{137.67}{11.00}{258}
\emline{145.33}{17.00}{259}{143.47}{14.85}{260}
\emline{143.47}{14.85}{261}{140.33}{13.34}{262}
\put(122.00,11.67){\vector(1,2){0.2}}
\emline{121.33}{10.33}{263}{122.00}{11.67}{264}
\put(150.33,17.00){\vector(1,1){0.2}}
\emline{149.67}{16.33}{265}{150.33}{17.00}{266}
\put(126.00,17.00){\makebox(0,0)[cc]{$=$}}
\put(51.00,18.00){\makebox(0,0)[cc]{,}}
\put(100.33,18.00){\makebox(0,0)[cc]{,}}
\put(152.67,15.00){\makebox(0,0)[cc]{.}}
\end{picture}
\\
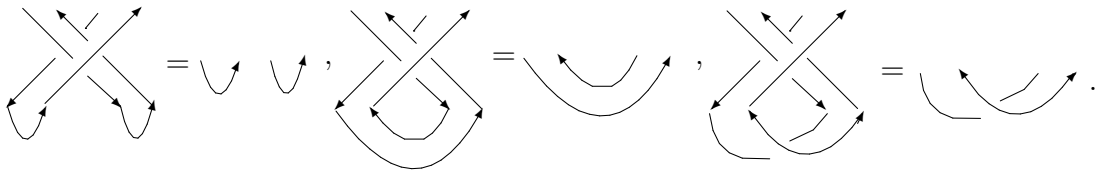
Figure\,2: $\sig$ acts as the identity on $\theta\ott\theta$, $\eta$, and $\xi$.
\end{center}
Consequently, $\gteni[2]\subseteq\ker A_2$ and $\gdsi[2]=\{0\}$. By
Corollary~\ref{c-1} $(\XX{\land}\XX)_0={0}$. 
\\
(iii) By \cite[(2.3)]{a-FRT},
$\Rda[wm]{vj}(\ell^{-\fettc})^m_{y}\ell^{+j}_z=\ell^{+w}_j(\ell^{-\fettc})^{v}_m\Rda[jy]{mz}$.
Further by \cite[Remark~21]{a-FRT},
\mbox{$\ell^\pm C^\fettt\ell^{\pm\fettt}=C^\fettt 1$.}
Let 
\begin{align*}
T_0=\ell^i_j\ell^m_nB^i_y\Rda[jy]{mz}C^n_z=(\ell^{-\fettc})^w_i\ell^{+w}_j(\ell^{-\fettc})^v_m\ell^{+v}_nB^i_y\Rda[jy]{mz}C^n_z.
\end{align*}
Using the above identities,  \rf[e-rc] twice, and finally \rf[e-x] one gets
\begin{align*}
T_0&=S\ell^{-i}_w(\ell^{-\fettc})^m_y\ell^{+j}_z\Rda[wm]{vj}\ell^{+v}_nB^i_yC^n_z
\\
&=S\ell^{-i}_wS\ell^{-y}_m\Rda[wm]{vj}C^v_jB^i_y
\\
&=\rtm S(\ell^{-y}_m\ell^{-i}_wC^w_m)B^i_y
\\
&=\rtm C^i_yB^i_y1=\rtm \xx.
\end{align*}
Using $X^\pm_{ij}=\ve_\pm\ell^i_j-\delta^i_j$, the above calculation, and again \rf[e-rc]
and \rf[e-x] 
it  follows that 
\begin{align*}
\mu(T)&=(\ve_\pm\ell^i_j-\delta^i_j)(\ve_\pm\ell^m_n-\delta^m_n)B^i_y\Rda[jy]{mz}C^n_z
\\
&=\rtm \xx-\ve_\pm B^i_y\Rda[iy]{mz}C^n_z\ell^m_n-\ve_\pm B^i_y\Rda[jy]{mz}C^m_z\ell^i_j+B^i_y\Rda[iy]{mz}C^m_z
\\
&=\rtm (\xx-\ve_\pm B^m_zC^n_z\ell^m_n-\ve_\pm B^i_yC^j_y\ell^i_j+\xx)
\\
&=-2\rtm X^\pm_0.
\end{align*}
Consequently, $T\in\mu^{-1}(\XX)$ and $\dim\mu^{-1}(\XX)_0\ge1$.
By Corollary\,\ref{c-1} applied to the projection
$P\in\Mor(\dr{\uhr}\gdul[2])$ onto the space $\gdui[2]$, the pairing
$$
\gdui[2]\times \mu^{-1}(\XX)_0\to \C
$$
is non-degenerate. Since $\gdui[2]=\gteni[2]{/}\SS(\RR)_\ii$,
$\dim\gteni[2]=3$ by (i), and $\dim \SS(\RR)_\ii\ge2$ by the result of Section\,\ref{uni}
we get $\dim\mu^{-1}(\XX)_0=\dim\gdui[2]=1$. This completes the proof.
\end{proof}

{\em Bicovariant subbimodules.}
We shall describe a method to construct a  class of bicovariant subbimodules 
of $\Gamm\ota\Gamm$. This method  is also applicable to higher tensor products 
$\gten[k]$, cf. \cite[p.\,1356]{a-HSu1}.
In this subsection $\AA$ is one of the Hopf algebras $\ooq$ or $\ospq$ and $\Gamm$ always denotes one of the $N^2$-dimensional 
bicovariant bimodules $\gpm$ over $\AA$ with left-invariant basis 
$\{\theta_{ij}|i,j=1,\dots,N\}$. The canonical basis of $\C^N$ is $\{e_1,\dots,e_N\}$.
\begin{lemma}             \label{l-sub}
{\em (i)} Let $\hat{P}$ and $\check{Q}$ be idempotents  in $\Mor(u\ott u)$ resp.
$\Mor(u^\fettc\ott u^\fettc)$.  We identify the underlying spaces of the right 
coaction on $\gl\ot\gl$ and the equivalent matrix corepresentation $u\ott 
u^\fettc \ott u\ott u^\fettc$ via $\theta_{ij}\ott\theta_{kl}\to e_i\ott 
e_j\ott e_k\ott e_l$. Then the subspace
\begin{align}\label{e-sub}
\Rlim{23}\hat{P}_{12}\check{Q}_{34}\Rre{23}(\gl\ott \gl)
\end{align}
of $\gl\ot\gl$ is the left-invariant basis of a 
bicovariant subbimodule of $\Gamm\ota\Gamm$ of dimension $\rank(\pda[])\rank(\check{Q})$.
\\
{\em (ii)} $\Gamm\ota\Gamm$ is the direct sum of $9$ bicovariant subbimodules
$\Lam[\tau\nu]$, $\tau,\nu\in\{+,-,0\}$, generated by the left-invariant
elements
\begin{align*}
\Rlim{23}\pda_{12}\pch[\nu]_{34}\Rre{23}(\gl\ott\gl).
\end{align*}
Moreover we have the following identity of bicovariant bimodules
$$
\ker A_2=\Lam[++]\oplus\Lam[--]\oplus\Lam[00].
$$
\end{lemma}
\begin{proof} (i) Since all four mappings appearing in \rf[e-sub] are morphisms 
of corepresentations one easily checks that 
$T:=\Rlim{23}\hat{P}_{12}\check{Q}_{34}\Rre{23}\in\Mor(u\ott u^\fettc\ott u\ott
u^\fettc)$. Hence the space is closed under the right coaction.
\\
Now we compute the right adjoint action. Set 
$\theta^{PQ}_{mnkl}=T(\theta_{mn}\ott\theta_{kl})$.
 By \rf[e-thetau] 
\begin{align*}
\theta^{PQ}_{mnkl}\rac u^i_j
&=T^{vwcd}_{mnkl}\Rda[ia]{vy}\Rda[wy]{eb}\Rda[es]{cz}\Rda[dz]{jt}
\theta_{ab}\ott\theta_{st}
\\
&=(\Rda{12}\Rre{23}\Rda{34}\Rre{45}T_{1234})^{iabst}_{mnklj}
\theta_{ab}\ott\theta_{st}
\\
&=(T_{1234}\Rda{12}\Rre{23}\Rda{34}\Rre{45})^{iabst}_{mnklj}\theta_{ab}\ott\theta_{st}
\\
&=\theta_{vwst}^{PQ}(\Rda{12}\Rre{23}\Rda{34}\Rre{45})^{ivwst}_{mnklj}.
\end{align*}
The  last  but one equation becomes evident  by taking a look at the graphical 
presentation of these equations.
\begin{center}
\special{em:linewidth 0.4pt}
\unitlength 1.00mm
\linethickness{0.4pt}
\begin{picture}(144.33,71.00)
\put(17.00,19.67){\framebox(15.00,10.00)[cc]{$\pda[]$}}
\put(38.67,19.67){\framebox(15.00,10.00)[cc]{$\check{Q}$}}
\put(19.00,9.33){\vector(0,-1){0.20}}
\emline{19.00}{19.33}{1}{19.00}{9.33}{2}
\put(51.67,19.67){\vector(0,1){0.20}}
\emline{51.67}{9.33}{3}{51.67}{19.67}{4}
\put(40.67,19.33){\vector(1,1){0.20}}
\emline{30.33}{9.00}{5}{40.67}{19.33}{6}
\emline{30.00}{19.67}{7}{34.67}{15.00}{8}
\put(40.33,9.67){\vector(1,-1){0.20}}
\emline{36.67}{13.33}{9}{40.33}{9.67}{10}
\emline{40.67}{39.00}{11}{36.67}{35.00}{12}
\put(30.33,29.67){\vector(-1,-1){0.20}}
\emline{34.33}{33.67}{13}{30.33}{29.67}{14}
\emline{39.33}{48.67}{15}{39.33}{48.67}{16}
\emline{28.67}{60.00}{17}{25.67}{57.00}{18}
\put(18.67,50.00){\vector(-1,-1){0.20}}
\emline{24.00}{55.33}{19}{18.67}{50.00}{20}
\emline{61.67}{37.33}{21}{61.67}{37.33}{22}
\emline{51.00}{49.00}{23}{48.00}{46.00}{24}
\put(41.00,39.00){\vector(-1,-1){0.20}}
\emline{46.33}{44.33}{25}{41.00}{39.00}{26}
\put(61.00,59.67){\vector(0,1){0.20}}
\emline{61.00}{48.67}{27}{61.00}{59.67}{28}
\emline{50.67}{59.67}{29}{50.67}{48.67}{30}
\put(19.00,30.00){\vector(0,-1){0.20}}
\emline{19.00}{49.67}{31}{19.00}{30.00}{32}
\put(30.33,40.00){\vector(-1,1){0.20}}
\emline{40.67}{29.67}{33}{30.33}{40.00}{34}
\emline{18.67}{60.00}{35}{30.00}{53.33}{36}
\emline{30.00}{40.33}{37}{30.00}{49.33}{38}
\put(39.33,59.67){\vector(1,1){0.20}}
\emline{30.00}{49.33}{39}{39.33}{59.67}{40}
\emline{51.67}{29.67}{41}{51.67}{38.67}{42}
\put(61.00,49.33){\vector(1,1){0.20}}
\emline{51.67}{39.00}{43}{61.00}{49.33}{44}
\emline{34.67}{51.33}{45}{52.00}{42.33}{46}
\emline{55.00}{40.33}{47}{61.00}{37.33}{48}
\put(60.67,9.00){\vector(0,-1){0.20}}
\emline{60.67}{37.33}{49}{60.67}{9.00}{50}
\emline{108.00}{20.33}{51}{108.00}{20.33}{52}
\emline{130.33}{9.33}{53}{130.33}{9.33}{54}
\emline{119.67}{20.67}{55}{116.67}{17.67}{56}
\put(109.67,10.67){\vector(-1,-1){0.20}}
\emline{115.00}{16.00}{57}{109.67}{10.67}{58}
\emline{119.33}{31.33}{59}{119.33}{20.33}{60}
\emline{87.33}{31.67}{61}{98.67}{25.00}{62}
\emline{100.33}{11.67}{63}{100.33}{20.67}{64}
\put(109.67,31.00){\vector(1,1){0.20}}
\emline{100.33}{20.67}{65}{109.67}{31.00}{66}
\emline{103.33}{23.00}{67}{120.67}{14.00}{68}
\put(96.67,42.00){\framebox(15.00,10.00)[cc]{$\pda[]$}}
\put(118.33,42.00){\framebox(15.00,10.00)[cc]{$\check{Q}$}}
\put(98.67,31.67){\vector(0,-1){0.20}}
\emline{98.67}{41.67}{69}{98.67}{31.67}{70}
\put(131.33,42.00){\vector(0,1){0.20}}
\emline{131.33}{31.67}{71}{131.33}{42.00}{72}
\put(120.33,41.67){\vector(1,1){0.20}}
\emline{110.00}{31.33}{73}{120.33}{41.67}{74}
\emline{109.67}{42.00}{75}{114.33}{37.33}{76}
\put(120.00,32.00){\vector(1,-1){0.20}}
\emline{116.33}{35.67}{77}{120.00}{32.00}{78}
\emline{120.33}{61.33}{79}{116.33}{57.33}{80}
\put(110.00,52.00){\vector(-1,-1){0.20}}
\emline{114.00}{56.00}{81}{110.00}{52.00}{82}
\emline{119.00}{71.00}{83}{119.00}{71.00}{84}
\emline{144.33}{59.67}{85}{144.33}{59.67}{86}
\put(110.00,62.33){\vector(-1,1){0.20}}
\emline{120.33}{52.00}{87}{110.00}{62.33}{88}
\put(99.00,52.00){\vector(0,-1){0.20}}
\emline{99.00}{62.33}{89}{99.00}{52.00}{90}
\put(87.67,31.67){\vector(0,-1){0.20}}
\emline{87.67}{62.33}{91}{87.67}{31.67}{92}
\emline{119.67}{10.33}{93}{131.33}{22.00}{94}
\emline{131.33}{22.00}{95}{131.33}{31.00}{96}
\emline{99.00}{31.33}{97}{95.67}{28.33}{98}
\emline{93.00}{26.33}{99}{87.67}{22.00}{100}
\put(87.67,10.00){\vector(0,-1){0.20}}
\emline{87.67}{22.00}{101}{87.67}{10.00}{102}
\put(129.67,10.00){\vector(3,-1){0.20}}
\emline{123.67}{12.33}{103}{129.67}{10.00}{104}
\put(72.67,37.00){\makebox(0,0)[cc]{$\theta_{ab}\ott\theta_{st}=$}}
\put(131.33,62.67){\vector(0,1){0.20}}
\emline{131.33}{52.00}{105}{131.33}{62.67}{106}
\put(18.67,63.00){\makebox(0,0)[cc]{$i$}}
\put(28.67,63.00){\makebox(0,0)[cc]{$a$}}
\put(39.67,63.00){\makebox(0,0)[cc]{$b$}}
\put(50.67,63.00){\makebox(0,0)[cc]{$s$}}
\put(61.00,63.00){\makebox(0,0)[cc]{$t$}}
\put(87.67,65.00){\makebox(0,0)[cc]{$i$}}
\put(99.00,65.00){\makebox(0,0)[cc]{$a$}}
\put(110.00,65.00){\makebox(0,0)[cc]{$b$}}
\put(120.00,65.00){\makebox(0,0)[cc]{$s$}}
\put(131.67,65.00){\makebox(0,0)[cc]{$t$}}
\put(19.00,6.00){\makebox(0,0)[cc]{$m$}}
\put(30.00,6.00){\makebox(0,0)[cc]{$n$}}
\put(40.67,6.00){\makebox(0,0)[cc]{$k$}}
\put(51.67,6.00){\makebox(0,0)[cc]{$l$}}
\put(60.67,6.00){\makebox(0,0)[cc]{$j$}}
\put(87.67,6.00){\makebox(0,0)[cc]{$m$}}
\put(100.33,6.00){\makebox(0,0)[cc]{$n$}}
\put(109.33,6.00){\makebox(0,0)[cc]{$k$}}
\put(119.67,6.00){\makebox(0,0)[cc]{$l$}}
\put(130.00,6.00){\makebox(0,0)[cc]{$j$}}
\put(142.33,37.33){\makebox(0,0)[cc]{$\theta_{ab}\ott\theta_{st}$.}}
\end{picture}
\\
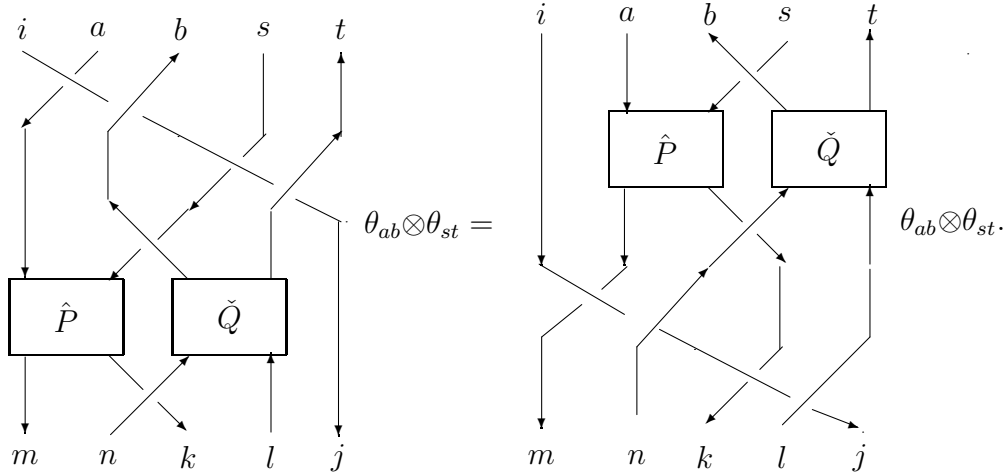
Figure\,3: $\langle\theta^{PQ}\rangle$ is closed under the right adjoint action.
\end{center}
Consequently $T(\gl\ott\gl)$ is closed under the right adjoint action. Hence
$\AA T(\gl\ott\gl)$ is a bicovariant subbimodule.
\\
(ii) The first part follows from (i) and the fact that
$(\pda[+]+\pda[-]+\pda[0])_{12}(\pch[+]+\pch[-]+\pch[0])_{34}$  is the
identity of $(\C^N)^{\ott4}$. In addition $\pda$ and $\pda[\nu]$ as well as
$\pch$ and $\pch[\nu]$ are pairwise orthogonal idempotents,
respectively. Hence the sum is direct. 
Let us turn to the second  part. Let $\lam_\tau$ denote the eigenvalue of $\Rda{}$ with
respect to the idempotent $\pda$, $\Rda{}\pda=\lamt\pda$, namely $\lam_+=q$, $\lam_-=-\qm$, and
$\lam_0=\rtm$. Note that $\Rch{}\pch=\lamt\pch$ as well. Put
$\rho=\Rlim{23}\pda_{12}\pch[\nu]_{34}\Rre{23}(\theta_{ijkl})$. 
Then by \rf[sig]
\begin{align*}
\sig(\rho)&= \Rlim{23}\Rda{12}\Rchm[34]\Rre{23}{\cdot}
\Rlim{23}\pda_{12}\pch[\nu]_{34}\Rre{23}(\theta_{ijkl})
\\
&= \lam_\tau\lam_\nu^{-1}\Rlim{23}\pda_{12}\pch[\nu]_{34}\Rre{23}(\theta_{ijkl})
\\
&= \lam_\tau\lam_\nu^{-1}\rho.
\end{align*}
Since $q$ is not a root of unity, $\lam_\tau\ne\lam_\nu$ for
$\tau\ne\nu$. Hence $\rho\in(\ker A_2)_\ll$ if and only if $\tau=\nu$. We
thus get $(\ker A_2)_\ll=\Lam[++]_\ll\oplus\Lam[--]_\ll\oplus
\Lam[00]_\ll$ as linear spaces.  By (i) each space on the right hand side generates a
bicovariant subbimodule. This completes the proof.
\end{proof}

To simplify notations we choose a new basis of $(\Gamm\ota\Gamm)_\ll$
\begin{align}\label{e-thetaq}
\thq_{vwst}=\Rlim[yz]{ws}\theta_{vy}\ott\theta_{zt},\qquad\theta_{vwst}=\Rre[yz]{ws}\thq_{vyzt}.
\end{align}
The right coaction now reads  $u\ott u\ott u^\fettc\ott u^\fettc$ and the
braiding in the new basis is $\ov{\sig}=\Rda{12}\Rchm[34]$. 
We simply write $\Lam[\tau]$ instead of $\Lam[\tau\tau]$, $\tau\in\{+,-,0\}$. 
Since the corresponding $\pda$ subcorepresentation of $u\ott u$ is
irreducible, by Schur's lemma $\Lam[\tau]$ has a {\em unique} up to scalars
bi-invariant element $\etat$. 
\begin{center}
\special{em:linewidth 0.4pt}
\unitlength 1.00mm
\linethickness{0.4pt}
\begin{picture}(105.00,31.00)
\put(25.00,15.00){\framebox(11.33,10.00)[cc]{$\pda$}}
\put(41.67,15.00){\framebox(11.33,10.00)[cc]{$\pch$}}
\put(27.00,25.00){\vector(0,-1){0.2}}
\emline{27.00}{28.33}{1}{27.00}{25.00}{2}
\put(33.67,25.00){\vector(0,-1){0.2}}
\emline{33.67}{28.00}{3}{33.67}{25.00}{4}
\put(44.00,27.67){\vector(0,1){0.2}}
\emline{44.00}{25.00}{5}{44.00}{27.67}{6}
\put(50.34,28.00){\vector(0,1){0.2}}
\emline{50.34}{25.00}{7}{50.34}{28.00}{8}
\put(44.00,15.00){\vector(3,4){0.2}}
\emline{33.67}{15.00}{9}{35.06}{13.25}{10}
\emline{35.06}{13.25}{11}{36.47}{12.06}{12}
\emline{36.47}{12.06}{13}{37.89}{11.44}{14}
\emline{37.89}{11.44}{15}{39.32}{11.38}{16}
\emline{39.32}{11.38}{17}{40.77}{11.89}{18}
\emline{40.77}{11.89}{19}{42.23}{12.96}{20}
\emline{42.23}{12.96}{21}{44.00}{15.00}{22}
\put(50.34,15.00){\vector(1,1){0.2}}
\emline{27.00}{15.00}{23}{28.89}{13.27}{24}
\emline{28.89}{13.27}{25}{30.76}{11.82}{26}
\emline{30.76}{11.82}{27}{32.60}{10.64}{28}
\emline{32.60}{10.64}{29}{34.41}{9.74}{30}
\emline{34.41}{9.74}{31}{36.20}{9.11}{32}
\emline{36.20}{9.11}{33}{37.95}{8.75}{34}
\emline{37.95}{8.75}{35}{39.67}{8.67}{36}
\emline{39.67}{8.67}{37}{41.37}{8.86}{38}
\emline{41.37}{8.86}{39}{43.04}{9.33}{40}
\emline{43.04}{9.33}{41}{44.68}{10.07}{42}
\emline{44.68}{10.07}{43}{46.29}{11.08}{44}
\emline{46.29}{11.08}{45}{47.87}{12.37}{46}
\emline{47.87}{12.37}{47}{50.34}{15.00}{48}
\put(3.33,20.33){\makebox(0,0)[cc]{$\etat=(\pda)^{mn}_{lk}\thq_{mnkl}=$}}
\put(61.33,20.33){\makebox(0,0)[cc]{$\thq_{mnkl}=$}}
\put(27.00,29.67){\makebox(0,0)[cc]{$m$}}
\put(33.67,29.67){\makebox(0,0)[cc]{$n$}}
\put(44.00,30.00){\makebox(0,0)[cc]{$k$}}
\put(50.34,30.00){\makebox(0,0)[cc]{$l$}}
\put(71.32,15.33){\framebox(11.33,10.00)[cc]{$\pda$}}
\put(73.33,25.33){\vector(0,-1){0.2}}
\emline{73.33}{28.67}{49}{73.33}{25.33}{50}
\put(79.99,25.33){\vector(0,-1){0.2}}
\emline{79.99}{28.33}{51}{79.99}{25.33}{52}
\put(90.33,15.33){\vector(3,4){0.2}}
\emline{79.99}{15.33}{53}{81.39}{13.58}{54}
\emline{81.39}{13.58}{55}{82.80}{12.39}{56}
\emline{82.80}{12.39}{57}{84.22}{11.77}{58}
\emline{84.22}{11.77}{59}{85.65}{11.71}{60}
\emline{85.65}{11.71}{61}{87.10}{12.22}{62}
\emline{87.10}{12.22}{63}{88.56}{13.29}{64}
\emline{88.56}{13.29}{65}{90.33}{15.33}{66}
\put(96.66,15.33){\vector(1,1){0.2}}
\emline{73.33}{15.33}{67}{75.23}{13.61}{68}
\emline{75.23}{13.61}{69}{77.10}{12.15}{70}
\emline{77.10}{12.15}{71}{78.94}{10.98}{72}
\emline{78.94}{10.98}{73}{80.75}{10.07}{74}
\emline{80.75}{10.07}{75}{82.53}{9.44}{76}
\emline{82.53}{9.44}{77}{84.28}{9.09}{78}
\emline{84.28}{9.09}{79}{86.01}{9.01}{80}
\emline{86.01}{9.01}{81}{87.70}{9.20}{82}
\emline{87.70}{9.20}{83}{89.37}{9.66}{84}
\emline{89.37}{9.66}{85}{91.01}{10.40}{86}
\emline{91.01}{10.40}{87}{92.62}{11.41}{88}
\emline{92.62}{11.41}{89}{94.20}{12.70}{90}
\emline{94.20}{12.70}{91}{96.66}{15.33}{92}
\put(73.33,30.00){\makebox(0,0)[cc]{$m$}}
\put(79.99,30.00){\makebox(0,0)[cc]{$n$}}
\put(90.33,28.33){\vector(0,1){0.2}}
\emline{90.33}{15.00}{93}{90.33}{28.33}{94}
\put(96.99,28.33){\vector(0,1){0.2}}
\emline{96.99}{15.00}{95}{96.99}{28.33}{96}
\put(105.00,20.67){\makebox(0,0)[cc]{$\thq_{mnkl}.$}}
\put(90.66,30.67){\makebox(0,0)[cc]{$k$}}
\put(97.33,31.00){\makebox(0,0)[cc]{$l$}}
\end{picture}
\\
Figure\,4: The bi-invariant elements $\eta^+$, $\eta^-$, and $\eta^0$.
\end{center}
The relations with the old basis of $\gteni[2]$ are
$\eta=\rt(\etat[0]+\etat[+]+\etat[-])$, $\xi=\xx\etat[0]$, and
\begin{align}\label{e-trans}
\theta\ott\theta=q\etat[+]-\qm\etat[-]+\rtm\etat[0].
\end{align}
The next lemma is the key step in our proof. 
\begin{lemma}\label{l-generator}
Let $\Gamm$ be one of the bicovariant FODC $\Gamm_\pm$ over $\AA$ and $\Lam[\tau]$, 
$\tau\in\{+,-,0\}$
the above defined bicovariant subbimodule of $\ker A_2$. 
Then $\Lam[\tau]$ is generated by the single element $\etat$. More precisely,
\begin{align*}
\etat\rac\AA=\Lam[\tau]_\ll.
\end{align*}
\end{lemma}
\begin{proof}
By Lemma~\ref{l-sub} the canonical left-invariant basis of $\Lam[\tau]$ is 
$$
\theta^\tau_{mnkl}=\Rlim{23}(\pda)_{12}(\pch)_{34}\Rre{23}
(\theta_{mn}\ott\theta_{kl})=\Rlim[ws]{yz}(\pda)^{vy}_{ma}
(\pch)^{zt}_{bl}\Rre[ab]{nk}\theta_{vw}\ott\theta_{st}.
$$
The proof is in two steps. First we compute $\etat\rac u^i_j$ and obtain
elements
\begin{align*}
\etat_{ij}&=B^i_z(\pda)^{mn}_{zk}(\pch)^{vw}_{ky}C^j_y\thq_{mnvw},
\\
\xi^\tau_{ij}&=(\pda)^{mn}_{yk}(\pch)^{vw}_{kz}\Rre[yz]{ij}\thq_{mnvw}.
\end{align*}
The graphical presentation of $\etat_{ij}$ and $\xi^\tau_{ij}$ is as follows.
\begin{center}
\special{em:linewidth 0.4pt}
\unitlength 1.00mm
\linethickness{0.4pt}
\begin{picture}(122.66,30.00)
\put(25.00,15.00){\framebox(11.33,10.00)[cc]{$\pda$}}
\put(41.67,15.00){\framebox(11.33,10.00)[cc]{$\pch$}}
\put(27.00,25.00){\vector(0,-1){0.2}}
\emline{27.00}{28.33}{1}{27.00}{25.00}{2}
\put(33.67,25.00){\vector(0,-1){0.2}}
\emline{33.67}{28.00}{3}{33.67}{25.00}{4}
\put(44.00,27.67){\vector(0,1){0.2}}
\emline{44.00}{25.00}{5}{44.00}{27.67}{6}
\put(50.34,28.00){\vector(0,1){0.2}}
\emline{50.34}{25.00}{7}{50.34}{28.00}{8}
\put(44.00,15.00){\vector(3,4){0.2}}
\emline{33.67}{15.00}{9}{35.07}{13.25}{10}
\emline{35.07}{13.25}{11}{36.47}{12.06}{12}
\emline{36.47}{12.06}{13}{37.89}{11.44}{14}
\emline{37.89}{11.44}{15}{39.33}{11.38}{16}
\emline{39.33}{11.38}{17}{40.77}{11.89}{18}
\emline{40.77}{11.89}{19}{42.23}{12.96}{20}
\emline{42.23}{12.96}{21}{44.00}{15.00}{22}
\put(17.66,20.33){\makebox(0,0)[cc]{$\etat_{ij}=$}}
\put(67.00,20.33){\makebox(0,0)[cc]{$\thq_{mnkl}$,\quad $\xi_{ij}^\tau=$}}
\put(27.00,29.67){\makebox(0,0)[cc]{$m$}}
\put(33.67,29.67){\makebox(0,0)[cc]{$n$}}
\put(44.00,30.00){\makebox(0,0)[cc]{$k$}}
\put(50.34,30.00){\makebox(0,0)[cc]{$l$}}
\put(27.00,12.33){\vector(0,-1){0.2}}
\emline{27.00}{15.00}{23}{27.00}{12.33}{24}
\put(27.00,10.33){\vector(0,1){0.2}}
\emline{27.00}{7.67}{25}{27.00}{10.33}{26}
\put(50.33,15.00){\vector(0,1){0.2}}
\emline{50.33}{11.33}{27}{50.33}{15.00}{28}
\put(50.33,7.67){\vector(0,-1){0.2}}
\emline{50.33}{11.33}{29}{50.33}{7.67}{30}
\put(27.00,4.67){\makebox(0,0)[cc]{$i$}}
\put(50.33,4.33){\makebox(0,0)[cc]{$j$}}
\put(86.33,15.00){\framebox(11.33,10.00)[cc]{$\pda$}}
\put(103.00,15.00){\framebox(11.33,10.00)[cc]{$\pch$}}
\put(88.33,25.00){\vector(0,-1){0.2}}
\emline{88.33}{28.33}{31}{88.33}{25.00}{32}
\put(95.00,25.00){\vector(0,-1){0.2}}
\emline{95.00}{28.00}{33}{95.00}{25.00}{34}
\put(105.33,27.67){\vector(0,1){0.2}}
\emline{105.33}{25.00}{35}{105.33}{27.67}{36}
\put(111.67,28.00){\vector(0,1){0.2}}
\emline{111.67}{25.00}{37}{111.67}{28.00}{38}
\put(105.33,15.00){\vector(3,4){0.2}}
\emline{95.00}{15.00}{39}{96.40}{13.25}{40}
\emline{96.40}{13.25}{41}{97.80}{12.06}{42}
\emline{97.80}{12.06}{43}{99.22}{11.44}{44}
\emline{99.22}{11.44}{45}{100.66}{11.38}{46}
\emline{100.66}{11.38}{47}{102.10}{11.89}{48}
\emline{102.10}{11.89}{49}{103.56}{12.96}{50}
\emline{103.56}{12.96}{51}{105.33}{15.00}{52}
\put(122.66,20.33){\makebox(0,0)[cc]{$\thq_{mnkl}$.}}
\put(88.33,29.67){\makebox(0,0)[cc]{$m$}}
\put(95.00,29.67){\makebox(0,0)[cc]{$n$}}
\put(105.33,30.00){\makebox(0,0)[cc]{$k$}}
\put(111.67,30.00){\makebox(0,0)[cc]{$l$}}
\put(88.33,4.67){\makebox(0,0)[cc]{$i$}}
\put(111.67,4.33){\makebox(0,0)[cc]{$j$}}
\emline{88.33}{7.33}{53}{97.67}{7.67}{54}
\emline{97.67}{7.67}{55}{101.33}{8.67}{56}
\emline{101.33}{8.67}{57}{104.67}{10.00}{58}
\emline{88.33}{15.00}{59}{98.33}{9.67}{60}
\put(111.67,7.00){\vector(1,0){0.2}}
\emline{102.33}{7.67}{61}{111.67}{7.00}{62}
\put(111.67,15.00){\vector(4,3){0.2}}
\emline{104.67}{10.00}{63}{111.67}{15.00}{64}
\end{picture}
\\
Figure\,5: The elements $\etat_{ij}$ and $\xi^\tau_{ij}$.
\end{center}
First we will show that
\begin{align}\label{e-etatu}
\etat\rac \plus{u^i_j}&=\QM(\lam_\tau^2+1)\xi^\tau_{ij}-
\rtm\QM(1+\lam_\tau^{-2})\etat_{ij},
\\
\etat\rac \plus{U}&=\alpt\etat,\quad\alpt=
\QM(\lam_\tau+\lam_\tau^{-1})(\lam_\tau\rt-\lam_\tau^{-1}\rtm).
\label{e-etatuu}
\end{align}
By \rf[e-thetaq] and \rf[e-thetau] one has 
\begin{align}
\thq_{mnkl}\rac u^i_j&=\Rlim[ab]{nk}(\theta_{ma}\rac u^i_c)(\theta_{bl}\rac
u^c_j)    \nn
\\
&=\Rlim[ab]{nk}\Rda[iv]{my}\Rda[ay]{cd}
\Rda[cp]{bz}\Rda[lz]{jt}\theta_{vd}\ott\theta_{pt}    \nn
\\
&=(\Rre{34}\Rda{12}\Rre{23}\Rda{34}\Rre{45}\Rlim{23})^{ivwst}_{mnklj}\thq_{vwst}.
\label{e-vwst}
\end{align}
\begin{center}
\special{em:linewidth 0.4pt}
\unitlength 1mm
\linethickness{0.4pt}
\begin{picture}(50.33,31.67)
\put(25.67,20.00){\vector(1,0){0.2}}
\emline{10.67}{10.00}{1}{11.36}{12.12}{2}
\emline{11.36}{12.12}{3}{12.28}{13.99}{4}
\emline{12.28}{13.99}{5}{13.43}{15.60}{6}
\emline{13.43}{15.60}{7}{14.81}{16.97}{8}
\emline{14.81}{16.97}{9}{16.42}{18.08}{10}
\emline{16.42}{18.08}{11}{18.26}{18.94}{12}
\emline{18.26}{18.94}{13}{20.32}{19.55}{14}
\emline{20.32}{19.55}{15}{22.61}{19.91}{16}
\emline{22.61}{19.91}{17}{25.67}{20.00}{18}
\put(16.67,10.00){\vector(0,-1){0.2}}
\emline{16.67}{16.67}{19}{16.67}{10.00}{20}
\put(21.67,10.00){\vector(0,-1){0.2}}
\emline{21.67}{16.67}{21}{21.67}{10.00}{22}
\put(28.67,28.00){\vector(0,1){0.2}}
\emline{28.67}{9.67}{23}{28.67}{28.00}{24}
\put(36.33,28.00){\vector(0,1){0.2}}
\emline{36.33}{10.00}{25}{36.33}{28.00}{26}
\put(45.33,10.00){\vector(-1,-4){0.2}}
\emline{39.00}{20.00}{27}{41.18}{19.14}{28}
\emline{41.18}{19.14}{29}{42.91}{18.02}{30}
\emline{42.91}{18.02}{31}{44.21}{16.64}{32}
\emline{44.21}{16.64}{33}{45.07}{15.00}{34}
\emline{45.07}{15.00}{35}{45.49}{13.10}{36}
\emline{45.49}{13.10}{37}{45.33}{10.00}{38}
\emline{21.67}{27.67}{39}{21.67}{21.33}{40}
\emline{21.67}{21.33}{41}{21.67}{21.33}{42}
\emline{21.67}{21.33}{43}{21.67}{21.33}{44}
\emline{21.67}{21.33}{45}{21.67}{21.33}{46}
\emline{16.67}{28.00}{47}{16.67}{20.00}{48}
\emline{30.67}{20.00}{49}{34.33}{20.00}{50}
\put(00.33,20.33){\makebox(0,0)[cc]{$\thq_{mnkl}\rac u^i_j=$}}
\put(10.67,6.33){\makebox(0,0)[cc]{$i$}}
\put(16.67,6.33){\makebox(0,0)[cc]{$m$}}
\put(21.67,6.33){\makebox(0,0)[cc]{$n$}}
\put(28.67,6.33){\makebox(0,0)[cc]{$k$}}
\put(36.33,6.33){\makebox(0,0)[cc]{$l$}}
\put(45.33,6.33){\makebox(0,0)[cc]{$j$}}
\put(16.67,31.00){\makebox(0,0)[cc]{$v$}}
\put(21.67,31.00){\makebox(0,0)[cc]{$w$}}
\put(28.33,31.00){\makebox(0,0)[cc]{$s$}}
\put(36.33,31.67){\makebox(0,0)[cc]{$t$}}
\put(50.33,20.33){\makebox(0,0)[cc]{$\thq_{vwst}$}}
\end{picture}
\\
Figure\,6: The right adjoint action of $u^i_j$ on $\thq_{mnkl}$.
\end{center}
In particular for $\etat$ we compute
\begin{center}
\special{em:linewidth 0.4pt}
\unitlength 1.00mm
\linethickness{0.4pt}
\begin{picture}(146.33,133.33)
\put(63.33,124.00){\framebox(11.67,7.00)[cc]{$\pda$}}
\put(79.00,124.00){\framebox(11.67,7.00)[cc]{$\pch$}}
\emline{60.66}{112.00}{1}{61.88}{114.00}{2}
\emline{61.88}{114.00}{3}{63.28}{115.74}{4}
\emline{63.28}{115.74}{5}{64.84}{117.24}{6}
\emline{64.84}{117.24}{7}{66.58}{118.49}{8}
\emline{66.58}{118.49}{9}{68.50}{119.48}{10}
\emline{68.50}{119.48}{11}{70.58}{120.23}{12}
\emline{70.58}{120.23}{13}{72.83}{120.73}{14}
\emline{72.83}{120.73}{15}{77.33}{121.00}{16}
\put(81.33,124.00){\vector(1,3){0.2}}
\emline{72.66}{119.33}{17}{74.85}{118.63}{18}
\emline{74.85}{118.63}{19}{76.73}{118.52}{20}
\emline{76.73}{118.52}{21}{78.33}{119.00}{22}
\emline{78.33}{119.00}{23}{79.62}{120.07}{24}
\emline{79.62}{120.07}{25}{80.63}{121.74}{26}
\emline{80.63}{121.74}{27}{81.33}{124.00}{28}
\put(88.33,123.67){\vector(1,3){0.2}}
\emline{66.00}{116.33}{29}{68.80}{115.60}{30}
\emline{68.80}{115.60}{31}{71.41}{115.07}{32}
\emline{71.41}{115.07}{33}{73.85}{114.74}{34}
\emline{73.85}{114.74}{35}{76.09}{114.61}{36}
\emline{76.09}{114.61}{37}{78.16}{114.69}{38}
\emline{78.16}{114.69}{39}{80.04}{114.97}{40}
\emline{80.04}{114.97}{41}{81.74}{115.45}{42}
\emline{81.74}{115.45}{43}{83.25}{116.13}{44}
\emline{83.25}{116.13}{45}{84.58}{117.01}{46}
\emline{84.58}{117.01}{47}{85.73}{118.10}{48}
\emline{85.73}{118.10}{49}{86.69}{119.38}{50}
\emline{86.69}{119.38}{51}{87.47}{120.87}{52}
\emline{87.47}{120.87}{53}{88.33}{123.67}{54}
\emline{66.00}{123.67}{55}{66.00}{119.33}{56}
\emline{72.00}{123.67}{57}{72.00}{121.67}{58}
\emline{81.33}{120.33}{59}{84.33}{119.00}{60}
\put(88.33,111.67){\vector(0,-1){0.2}}
\emline{86.66}{117.33}{61}{87.79}{116.31}{62}
\emline{87.79}{116.31}{63}{88.33}{111.67}{64}
\put(60.66,109.00){\makebox(0,0)[cc]{$i$}}
\put(88.33,108.67){\makebox(0,0)[cc]{$j$}}
\put(81.33,133.33){\vector(0,1){0.2}}
\emline{81.33}{131.00}{65}{81.33}{133.33}{66}
\put(87.66,133.33){\vector(0,1){0.2}}
\emline{87.66}{131.00}{67}{87.66}{133.33}{68}
\put(66.00,131.00){\vector(0,-1){0.2}}
\emline{66.00}{133.33}{69}{66.00}{131.00}{70}
\put(72.00,131.00){\vector(0,-1){0.2}}
\emline{72.00}{133.00}{71}{72.00}{131.00}{72}
\put(70.66,118.00){\dashbox{0.33}(3.67,4.67)[cc]{}}
\put(35.00,121.67){\makebox(0,0)[cc]
{$\etat\rac u^i_j=(\pda)^{mn}_{lk}\thq_{mnkl}\rac u^i_j=$}}
\put(97.00,122.00){\makebox(0,0)[cc]{$=$}}
\put(13.00,91.00){\framebox(11.67,7.00)[cc]{$\pda$}}
\put(28.67,91.00){\framebox(11.67,7.00)[cc]{$\pch$}}
\put(38.00,90.67){\vector(1,3){0.2}}
\emline{15.67}{83.33}{73}{18.47}{82.60}{74}
\emline{18.47}{82.60}{75}{21.08}{82.07}{76}
\emline{21.08}{82.07}{77}{23.52}{81.74}{78}
\emline{23.52}{81.74}{79}{25.76}{81.61}{80}
\emline{25.76}{81.61}{81}{27.83}{81.69}{82}
\emline{27.83}{81.69}{83}{29.71}{81.97}{84}
\emline{29.71}{81.97}{85}{31.41}{82.45}{86}
\emline{31.41}{82.45}{87}{32.92}{83.13}{88}
\emline{32.92}{83.13}{89}{34.25}{84.01}{90}
\emline{34.25}{84.01}{91}{35.40}{85.10}{92}
\emline{35.40}{85.10}{93}{36.36}{86.38}{94}
\emline{36.36}{86.38}{95}{37.14}{87.87}{96}
\emline{37.14}{87.87}{97}{38.00}{90.67}{98}
\emline{15.67}{90.67}{99}{15.67}{86.33}{100}
\put(38.00,78.67){\vector(0,-1){0.2}}
\emline{36.33}{84.33}{101}{37.46}{83.31}{102}
\emline{37.46}{83.31}{103}{38.00}{78.67}{104}
\put(10.33,76.00){\makebox(0,0)[cc]{$i$}}
\put(38.00,75.67){\makebox(0,0)[cc]{$j$}}
\put(31.00,100.33){\vector(0,1){0.2}}
\emline{31.00}{98.00}{105}{31.00}{100.33}{106}
\put(37.33,100.33){\vector(0,1){0.2}}
\emline{37.33}{98.00}{107}{37.33}{100.33}{108}
\put(15.67,98.00){\vector(0,-1){0.2}}
\emline{15.67}{100.33}{109}{15.67}{98.00}{110}
\put(21.67,98.00){\vector(0,-1){0.2}}
\emline{21.67}{100.00}{111}{21.67}{98.00}{112}
\put(57.66,91.00){\framebox(11.67,7.00)[cc]{$\pda$}}
\put(73.33,91.00){\framebox(11.67,7.00)[cc]{$\pch$}}
\put(82.66,90.67){\vector(1,3){0.2}}
\emline{60.33}{83.33}{113}{63.13}{82.60}{114}
\emline{63.13}{82.60}{115}{65.74}{82.07}{116}
\emline{65.74}{82.07}{117}{68.18}{81.74}{118}
\emline{68.18}{81.74}{119}{70.42}{81.61}{120}
\emline{70.42}{81.61}{121}{72.49}{81.69}{122}
\emline{72.49}{81.69}{123}{74.37}{81.97}{124}
\emline{74.37}{81.97}{125}{76.07}{82.45}{126}
\emline{76.07}{82.45}{127}{77.58}{83.13}{128}
\emline{77.58}{83.13}{129}{78.91}{84.01}{130}
\emline{78.91}{84.01}{131}{80.06}{85.10}{132}
\emline{80.06}{85.10}{133}{81.02}{86.38}{134}
\emline{81.02}{86.38}{135}{81.80}{87.87}{136}
\emline{81.80}{87.87}{137}{82.66}{90.67}{138}
\emline{60.33}{90.67}{139}{60.33}{86.33}{140}
\emline{75.66}{87.33}{141}{78.66}{86.00}{142}
\put(82.66,78.67){\vector(0,-1){0.2}}
\emline{80.99}{84.33}{143}{82.12}{83.31}{144}
\emline{82.12}{83.31}{145}{82.66}{78.67}{146}
\put(54.99,76.00){\makebox(0,0)[cc]{$i$}}
\put(82.66,75.67){\makebox(0,0)[cc]{$j$}}
\put(75.66,100.33){\vector(0,1){0.2}}
\emline{75.66}{98.00}{147}{75.66}{100.33}{148}
\put(81.99,100.33){\vector(0,1){0.2}}
\emline{81.99}{98.00}{149}{81.99}{100.33}{150}
\put(60.33,98.00){\vector(0,-1){0.2}}
\emline{60.33}{100.33}{151}{60.33}{98.00}{152}
\put(66.33,98.00){\vector(0,-1){0.2}}
\emline{66.33}{100.00}{153}{66.33}{98.00}{154}
\put(114.00,91.00){\framebox(11.67,7.00)[cc]{$\pda$}}
\put(129.67,91.00){\framebox(11.67,7.00)[cc]{$\pch$}}
\put(139.00,90.67){\vector(1,3){0.2}}
\emline{116.67}{83.33}{155}{119.47}{82.60}{156}
\emline{119.47}{82.60}{157}{122.08}{82.07}{158}
\emline{122.08}{82.07}{159}{124.52}{81.74}{160}
\emline{124.52}{81.74}{161}{126.76}{81.61}{162}
\emline{126.76}{81.61}{163}{128.83}{81.69}{164}
\emline{128.83}{81.69}{165}{130.71}{81.97}{166}
\emline{130.71}{81.97}{167}{132.41}{82.45}{168}
\emline{132.41}{82.45}{169}{133.92}{83.13}{170}
\emline{133.92}{83.13}{171}{135.25}{84.01}{172}
\emline{135.25}{84.01}{173}{136.40}{85.10}{174}
\emline{136.40}{85.10}{175}{137.36}{86.38}{176}
\emline{137.36}{86.38}{177}{138.14}{87.87}{178}
\emline{138.14}{87.87}{179}{139.00}{90.67}{180}
\emline{116.67}{90.67}{181}{116.67}{86.33}{182}
\put(139.00,78.67){\vector(0,-1){0.2}}
\emline{137.33}{84.33}{183}{138.46}{83.31}{184}
\emline{138.46}{83.31}{185}{139.00}{78.67}{186}
\put(111.33,76.00){\makebox(0,0)[cc]{$i$}}
\put(139.00,75.67){\makebox(0,0)[cc]{$j$}}
\put(132.00,100.33){\vector(0,1){0.2}}
\emline{132.00}{98.00}{187}{132.00}{100.33}{188}
\put(138.33,100.33){\vector(0,1){0.2}}
\emline{138.33}{98.00}{189}{138.33}{100.33}{190}
\put(116.67,98.00){\vector(0,-1){0.2}}
\emline{116.67}{100.33}{191}{116.67}{98.00}{192}
\put(122.67,98.00){\vector(0,-1){0.2}}
\emline{122.67}{100.00}{193}{122.67}{98.00}{194}
\put(5.00,89.33){\makebox(0,0)[cc]{$=$}}
\put(47.67,89.33){\makebox(0,0)[cc]{$+\QM$}}
\put(99.00,89.33){\makebox(0,0)[cc]{$-\QM$}}
\put(31.67,91.00){\vector(4,3){0.2}}
\emline{21.67}{91.00}{195}{23.23}{89.42}{196}
\emline{23.23}{89.42}{197}{24.80}{88.42}{198}
\emline{24.80}{88.42}{199}{26.36}{88.01}{200}
\emline{26.36}{88.01}{201}{27.92}{88.19}{202}
\emline{27.92}{88.19}{203}{29.48}{88.95}{204}
\emline{29.48}{88.95}{205}{31.67}{91.00}{206}
\emline{10.33}{79.00}{207}{11.56}{80.68}{208}
\emline{11.56}{80.68}{209}{12.93}{82.15}{210}
\emline{12.93}{82.15}{211}{14.43}{83.42}{212}
\emline{14.43}{83.42}{213}{16.07}{84.48}{214}
\emline{16.07}{84.48}{215}{17.84}{85.34}{216}
\emline{17.84}{85.34}{217}{19.75}{86.00}{218}
\emline{19.75}{86.00}{219}{21.79}{86.46}{220}
\emline{21.79}{86.46}{221}{23.96}{86.71}{222}
\emline{23.96}{86.71}{223}{26.27}{86.75}{224}
\emline{26.27}{86.75}{225}{28.71}{86.60}{226}
\emline{28.71}{86.60}{227}{31.29}{86.23}{228}
\emline{31.29}{86.23}{229}{34.00}{85.67}{230}
\emline{55.00}{79.33}{231}{56.09}{81.12}{232}
\emline{56.09}{81.12}{233}{57.42}{82.64}{234}
\emline{57.42}{82.64}{235}{58.99}{83.88}{236}
\emline{58.99}{83.88}{237}{60.80}{84.84}{238}
\emline{60.80}{84.84}{239}{62.85}{85.52}{240}
\emline{62.85}{85.52}{241}{65.15}{85.93}{242}
\emline{65.15}{85.93}{243}{69.33}{86.00}{244}
\put(75.67,91.00){\vector(1,2){0.2}}
\emline{68.67}{86.00}{245}{71.11}{86.56}{246}
\emline{71.11}{86.56}{247}{73.09}{87.59}{248}
\emline{73.09}{87.59}{249}{74.61}{89.06}{250}
\emline{74.61}{89.06}{251}{75.67}{91.00}{252}
\emline{66.33}{91.00}{253}{67.97}{89.39}{254}
\emline{67.97}{89.39}{255}{72.00}{88.00}{256}
\emline{132.33}{90.67}{257}{132.34}{88.42}{258}
\emline{132.34}{88.42}{259}{131.97}{86.68}{260}
\emline{131.97}{86.68}{261}{131.20}{85.44}{262}
\emline{131.20}{85.44}{263}{130.04}{84.71}{264}
\emline{130.04}{84.71}{265}{128.49}{84.47}{266}
\emline{128.49}{84.47}{267}{125.67}{85.00}{268}
\emline{127.00}{84.67}{269}{126.88}{86.48}{270}
\emline{126.88}{86.48}{271}{127.64}{87.47}{272}
\emline{127.64}{87.47}{273}{130.67}{87.33}{274}
\emline{133.00}{87.00}{275}{134.67}{86.00}{276}
\put(146.33,89.67){\makebox(0,0)[cc]{$=$}}
\put(13.33,80.67){\dashbox{0.67}(4.33,7.33)[cc]{}}
\put(57.67,81.00){\dashbox{0.67}(5.33,6.67)[cc]{}}
\put(13.00,58.33){\framebox(11.67,7.00)[cc]{$\pda$}}
\put(28.67,58.33){\framebox(11.67,7.00)[cc]{$\pch$}}
\put(10.33,43.33){\makebox(0,0)[cc]{$i$}}
\put(38.00,43.00){\makebox(0,0)[cc]{$j$}}
\put(31.00,67.66){\vector(0,1){0.2}}
\emline{31.00}{65.33}{277}{31.00}{67.66}{278}
\put(37.33,67.66){\vector(0,1){0.2}}
\emline{37.33}{65.33}{279}{37.33}{67.66}{280}
\put(15.67,65.33){\vector(0,-1){0.2}}
\emline{15.67}{67.66}{281}{15.67}{65.33}{282}
\put(21.67,65.33){\vector(0,-1){0.2}}
\emline{21.67}{67.33}{283}{21.67}{65.33}{284}
\put(5.00,56.66){\makebox(0,0)[cc]{$=$}}
\put(31.67,58.33){\vector(4,3){0.2}}
\emline{21.67}{58.33}{285}{23.23}{56.75}{286}
\emline{23.23}{56.75}{287}{24.80}{55.75}{288}
\emline{24.80}{55.75}{289}{26.36}{55.34}{290}
\emline{26.36}{55.34}{291}{27.92}{55.52}{292}
\emline{27.92}{55.52}{293}{29.48}{56.28}{294}
\emline{29.48}{56.28}{295}{31.67}{58.33}{296}
\put(37.33,58.33){\vector(3,4){0.2}}
\emline{15.67}{58.33}{297}{17.47}{56.56}{298}
\emline{17.47}{56.56}{299}{19.25}{55.07}{300}
\emline{19.25}{55.07}{301}{21.00}{53.88}{302}
\emline{21.00}{53.88}{303}{22.73}{52.98}{304}
\emline{22.73}{52.98}{305}{24.43}{52.37}{306}
\emline{24.43}{52.37}{307}{26.10}{52.05}{308}
\emline{26.10}{52.05}{309}{27.74}{52.02}{310}
\emline{27.74}{52.02}{311}{29.36}{52.28}{312}
\emline{29.36}{52.28}{313}{30.95}{52.84}{314}
\emline{30.95}{52.84}{315}{32.51}{53.68}{316}
\emline{32.51}{53.68}{317}{34.05}{54.81}{318}
\emline{34.05}{54.81}{319}{35.55}{56.24}{320}
\emline{35.55}{56.24}{321}{37.33}{58.33}{322}
\put(37.67,45.67){\vector(2,-1){0.2}}
\emline{10.33}{46.00}{323}{12.81}{46.91}{324}
\emline{12.81}{46.91}{325}{15.26}{47.66}{326}
\emline{15.26}{47.66}{327}{17.66}{48.23}{328}
\emline{17.66}{48.23}{329}{20.03}{48.63}{330}
\emline{20.03}{48.63}{331}{22.37}{48.86}{332}
\emline{22.37}{48.86}{333}{24.67}{48.92}{334}
\emline{24.67}{48.92}{335}{26.92}{48.80}{336}
\emline{26.92}{48.80}{337}{29.15}{48.52}{338}
\emline{29.15}{48.52}{339}{31.33}{48.06}{340}
\emline{31.33}{48.06}{341}{33.48}{47.44}{342}
\emline{33.48}{47.44}{343}{35.59}{46.64}{344}
\emline{35.59}{46.64}{345}{37.67}{45.67}{346}
\put(108.33,58.33){\framebox(11.67,7.00)[cc]{$\pda$}}
\put(124.00,58.33){\framebox(11.67,7.00)[cc]{$\pch$}}
\put(126.33,67.66){\vector(0,1){0.2}}
\emline{126.33}{65.33}{347}{126.33}{67.66}{348}
\put(132.66,67.66){\vector(0,1){0.2}}
\emline{132.66}{65.33}{349}{132.66}{67.66}{350}
\put(111.00,65.33){\vector(0,-1){0.2}}
\emline{111.00}{67.66}{351}{111.00}{65.33}{352}
\put(117.00,65.33){\vector(0,-1){0.2}}
\emline{117.00}{67.33}{353}{117.00}{65.33}{354}
\put(127.00,58.33){\vector(4,3){0.2}}
\emline{117.00}{58.33}{355}{118.56}{56.75}{356}
\emline{118.56}{56.75}{357}{120.13}{55.75}{358}
\emline{120.13}{55.75}{359}{121.69}{55.34}{360}
\emline{121.69}{55.34}{361}{123.25}{55.52}{362}
\emline{123.25}{55.52}{363}{124.81}{56.28}{364}
\emline{124.81}{56.28}{365}{127.00}{58.33}{366}
\put(110.66,55.00){\vector(0,-1){0.2}}
\emline{110.66}{58.33}{367}{110.66}{55.00}{368}
\put(110.66,53.00){\vector(0,1){0.2}}
\emline{110.66}{46.33}{369}{110.66}{53.00}{370}
\put(132.66,47.00){\vector(0,-1){0.2}}
\emline{132.66}{54.00}{371}{132.66}{47.00}{372}
\put(132.66,58.00){\vector(0,1){0.2}}
\emline{132.66}{53.67}{373}{132.66}{58.00}{374}
\put(98.66,58.00){\makebox(0,0)[cc]{$-\QM\rtm$}}
\put(58.33,58.33){\framebox(11.67,7.00)[cc]{$\pda$}}
\put(74.00,58.33){\framebox(11.67,7.00)[cc]{$\pch$}}
\put(76.33,67.66){\vector(0,1){0.2}}
\emline{76.33}{65.33}{375}{76.33}{67.66}{376}
\put(82.66,67.66){\vector(0,1){0.2}}
\emline{82.66}{65.33}{377}{82.66}{67.66}{378}
\put(61.00,65.33){\vector(0,-1){0.2}}
\emline{61.00}{67.66}{379}{61.00}{65.33}{380}
\put(67.00,65.33){\vector(0,-1){0.2}}
\emline{67.00}{67.33}{381}{67.00}{65.33}{382}
\put(77.00,58.33){\vector(4,3){0.2}}
\emline{67.00}{58.33}{383}{68.56}{56.75}{384}
\emline{68.56}{56.75}{385}{70.13}{55.75}{386}
\emline{70.13}{55.75}{387}{71.69}{55.34}{388}
\emline{71.69}{55.34}{389}{73.25}{55.52}{390}
\emline{73.25}{55.52}{391}{74.81}{56.28}{392}
\emline{74.81}{56.28}{393}{77.00}{58.33}{394}
\put(82.33,58.00){\vector(2,1){0.2}}
\emline{61.00}{47.00}{395}{82.33}{58.00}{396}
\emline{61.00}{58.33}{397}{70.00}{53.67}{398}
\put(83.00,47.33){\vector(2,-1){0.2}}
\emline{74.00}{51.67}{399}{83.00}{47.33}{400}
\put(110.66,43.00){\makebox(0,0)[cc]{$i$}}
\put(132.33,43.33){\makebox(0,0)[cc]{$j$}}
\put(60.67,43.67){\makebox(0,0)[cc]{$i$}}
\put(83.00,43.67){\makebox(0,0)[cc]{$j$}}
\put(48.33,58.33){\makebox(0,0)[cc]{$+\QM$}}
\put(146.00,58.33){\makebox(0,0)[cc]{$+$}}
\put(29.00,25.00){\framebox(11.67,7.00)[cc]{$\pda$}}
\put(44.67,25.00){\framebox(11.67,7.00)[cc]{$\pch$}}
\put(47.00,34.33){\vector(0,1){0.2}}
\emline{47.00}{32.00}{401}{47.00}{34.33}{402}
\put(53.33,34.33){\vector(0,1){0.2}}
\emline{53.33}{32.00}{403}{53.33}{34.33}{404}
\put(31.67,32.00){\vector(0,-1){0.2}}
\emline{31.67}{34.33}{405}{31.67}{32.00}{406}
\put(37.67,32.00){\vector(0,-1){0.2}}
\emline{37.67}{34.00}{407}{37.67}{32.00}{408}
\put(47.67,25.00){\vector(4,3){0.2}}
\emline{37.67}{25.00}{409}{39.23}{23.42}{410}
\emline{39.23}{23.42}{411}{40.79}{22.42}{412}
\emline{40.79}{22.42}{413}{42.35}{22.01}{414}
\emline{42.35}{22.01}{415}{43.92}{22.19}{416}
\emline{43.92}{22.19}{417}{45.48}{22.95}{418}
\emline{45.48}{22.95}{419}{47.67}{25.00}{420}
\put(53.00,24.67){\vector(2,1){0.2}}
\emline{31.67}{13.67}{421}{53.00}{24.67}{422}
\emline{31.67}{25.00}{423}{40.67}{20.33}{424}
\put(53.67,14.00){\vector(2,-1){0.2}}
\emline{44.67}{18.33}{425}{53.67}{14.00}{426}
\put(31.33,10.33){\makebox(0,0)[cc]{$i$}}
\put(53.67,10.33){\makebox(0,0)[cc]{$j$}}
\put(13.67,25.00){\makebox(0,0)[cc]{$+\lam_\tau^2\QM$}}
\put(80.67,25.00){\framebox(11.67,7.00)[cc]{$\pda$}}
\put(96.34,25.00){\framebox(11.67,7.00)[cc]{$\pch$}}
\put(98.67,34.33){\vector(0,1){0.2}}
\emline{98.67}{32.00}{427}{98.67}{34.33}{428}
\put(105.00,34.33){\vector(0,1){0.2}}
\emline{105.00}{32.00}{429}{105.00}{34.33}{430}
\put(83.34,32.00){\vector(0,-1){0.2}}
\emline{83.34}{34.33}{431}{83.34}{32.00}{432}
\put(89.34,32.00){\vector(0,-1){0.2}}
\emline{89.34}{34.00}{433}{89.34}{32.00}{434}
\put(99.34,25.00){\vector(4,3){0.2}}
\emline{89.34}{25.00}{435}{90.90}{23.42}{436}
\emline{90.90}{23.42}{437}{92.46}{22.42}{438}
\emline{92.46}{22.42}{439}{94.02}{22.01}{440}
\emline{94.02}{22.01}{441}{95.59}{22.19}{442}
\emline{95.59}{22.19}{443}{97.15}{22.95}{444}
\emline{97.15}{22.95}{445}{99.34}{25.00}{446}
\put(83.00,21.67){\vector(0,-1){0.2}}
\emline{83.00}{25.00}{447}{83.00}{21.67}{448}
\put(83.00,19.67){\vector(0,1){0.2}}
\emline{83.00}{13.00}{449}{83.00}{19.67}{450}
\put(105.00,13.67){\vector(0,-1){0.2}}
\emline{105.00}{20.67}{451}{105.00}{13.67}{452}
\put(105.00,24.67){\vector(0,1){0.2}}
\emline{105.00}{20.33}{453}{105.00}{24.67}{454}
\put(71.00,24.67){\makebox(0,0)[cc]{$-\QM\rtm\lam_\tau^{-2}$}}
\put(83.00,9.67){\makebox(0,0)[cc]{$i$}}
\put(104.67,10.00){\makebox(0,0)[cc]{$j$}}
\put(113.34,25.67){\makebox(0,0)[cc]{.}}
\put(119.33,85.00){\vector(-4,-3){0.2}}
\emline{122.67}{91.00}{455}{122.48}{88.55}{456}
\emline{122.48}{88.55}{457}{121.52}{86.63}{458}
\emline{121.52}{86.63}{459}{119.33}{85.00}{460}
\put(118.33,84.67){\vector(3,1){0.2}}
\emline{111.33}{78.67}{461}{112.16}{80.73}{462}
\emline{112.16}{80.73}{463}{113.44}{82.39}{464}
\emline{113.44}{82.39}{465}{115.17}{83.63}{466}
\emline{115.17}{83.63}{467}{118.33}{84.67}{468}
\end{picture}
\\
Figure\,7: The proof of \rf[e-etatu].
\end{center}
In the first step we replaced the crossing  in the dash box using \rf[e-rm]. In 
the second step we did the same with the $\Rda{}$-matrix in the first dash box. 
This gives the first three terms in the next line. 
Moreover the dash box in the second summand is multiplied by  $\pch$ and gives
$\lam_\tau I$ (no crossing). Similarly, a second crossing in the same term 
gives another $\lam_\tau$. With the  third summand we proceed in the same 
way; in addition the curl gives the factor $\rtm$. Since
$\plus{u^i_j}=u^i_j-\delta_{ij}1$, \rf[e-etatu] follows immediately. 
\\
Note that for $\tau=0$,
\begin{align}\label{e-tn}
\etat[0]\rac u^i_j=\delta_{ij}\etat[0]
\end{align}
is obvious from the first line in 
Figure\,7 since $\pda[0]=\xx^{-1}K$ and no
crossing appears there. Moreover $\etat[0]_{ij}=\delta_{ij}\etat[0]$ and
$\xi^0_{ij}=\rt\delta_{ij}\etat[0]$ and \rf[e-etatu] and \rf[e-etatuu] are
valid. Since $\Lam[0]$ is one-dimensional there is nothing to prove. Now we
fix $\tau\in\{+,-\}$. We shall eliminate $\etat_{ij}$ from \rf[e-etatu].
Multiplying \rf[e-etatu] by $D^j_i$ and 
using $D^j_i\etat_{ij}=\etat$, $D^j_i\xi^\tau_{ij}=\rt\etat$ gives
 \rf[e-etatuu]. 
Since $q$ is not a root of unity,
$T_\tau=\QM(\lam_\tau^2+1)\Rda[-1]{}-\rtm\QM(1+\lam_\tau^{-2})$ is invertible with inverse 
\begin{align*}
T_+^{-1}&=\frac{1}{\QM\qtwo}\bigl(\frac{1}{1-\qm\rtm}\pda[+]+\frac{1}{-q^2-\qm\rtm}\pda[-]+\frac{1}{\rt
  q-\rtm\qm}\pda[0]\bigr),
\\
T_-^{-1}&=\frac{1}{\QM\qtwo}\bigl(\frac{1}{q^{-2}-q\rtm}\pda[+]+\frac{1}{-1-q\rtm}\pda[-]+\frac{1}{\qm\rt-q\rtm}\pda[0]\bigr).
\end{align*}
 Set
$(S_\tau)^{ij}_{st}=B^y_i(T_\tau^{-1})^{yj}_{zt}C^s_z$ and multiply
\rf[e-etatu] by  $(S_\tau)^{ij}_{st}$. Then we obtain 
$(S_\tau)^{ij}_{st}\etat\rac\plus{u^i_j}=\etat_{st}$. 
Consequently, $\etat_{st}\in \Lam[\tau]$ for  $s,t=1,\dots,N$ and  
$\tau\in\{+,-\}$.

In the second step we again compute the right adjoint action of $u^i_j$ but on 
elements $\etat_{st}$. We obtain elements 
\begin{align*}
\xi^\tau_{sijt}&=(\pda)^{mn}_{sy}(\pch)^{vw}_{zt}\Rre[yz]{ij}\thq_{mnvw},
\\
\etat_{sijt}&=B^s_yB^i_z(\pda)^{mn}_{yz}(\pch)^{vw}_{dc}C^j_dC^t_c\thq_{mnvw}.
\end{align*}
Obviously $\xi^\tau_{sijt}=\etat_{syzt}\Rre[yz]{ij}$. Graphically they are 
represented by
\begin{center}
\special{em:linewidth 0.4pt}
\unitlength 1.00mm
\linethickness{0.4pt}
\begin{picture}(108.67,26.99)
\put(23.99,17.66){\framebox(11.67,7.00)[cc]{$\pda$}}
\put(39.66,17.66){\framebox(11.67,7.00)[cc]{$\pch$}}
\put(41.99,26.99){\vector(0,1){0.2}}
\emline{41.99}{24.99}{1}{41.99}{26.99}{2}
\put(48.32,26.99){\vector(0,1){0.2}}
\emline{48.32}{24.99}{3}{48.32}{26.99}{4}
\put(26.66,24.99){\vector(0,-1){0.2}}
\emline{26.66}{26.99}{5}{26.66}{24.99}{6}
\put(32.66,24.99){\vector(0,-1){0.2}}
\emline{32.66}{26.66}{7}{32.66}{24.99}{8}
\put(26.32,14.33){\vector(0,-1){0.2}}
\emline{26.32}{17.66}{9}{26.32}{14.33}{10}
\put(48.32,17.33){\vector(0,1){0.2}}
\emline{48.32}{13.00}{11}{48.32}{17.33}{12}
\put(26.32,5.33){\makebox(0,0)[cc]{$i$}}
\put(48.66,5.33){\makebox(0,0)[cc]{$j$}}
\put(41.99,17.00){\vector(4,3){0.2}}
\emline{32.66}{9.33}{13}{41.99}{17.00}{14}
\emline{32.66}{17.33}{15}{36.99}{14.00}{16}
\put(42.66,9.66){\vector(4,-3){0.2}}
\emline{37.99}{13.33}{17}{42.66}{9.66}{18}
\put(26.33,13.33){\vector(0,1){0.2}}
\emline{26.33}{9.66}{19}{26.33}{13.33}{20}
\put(32.66,5.33){\makebox(0,0)[cc]{$s$}}
\put(42.33,5.33){\makebox(0,0)[cc]{$t$}}
\put(17.66,17.66){\makebox(0,0)[cc]{$\xi^\tau_{istj}=$}}
\put(95.01,26.99){\vector(0,1){0.2}}
\emline{95.01}{24.99}{21}{95.01}{26.99}{22}
\put(101.34,26.99){\vector(0,1){0.2}}
\emline{101.34}{24.99}{23}{101.34}{26.99}{24}
\put(79.68,24.99){\vector(0,-1){0.2}}
\emline{79.68}{26.99}{25}{79.68}{24.99}{26}
\put(85.68,24.99){\vector(0,-1){0.2}}
\emline{85.68}{26.66}{27}{85.68}{24.99}{28}
\put(79.34,14.33){\vector(0,-1){0.2}}
\emline{79.34}{17.66}{29}{79.34}{14.33}{30}
\put(101.34,17.33){\vector(0,1){0.2}}
\emline{101.34}{13.00}{31}{101.34}{17.33}{32}
\put(79.34,5.33){\makebox(0,0)[cc]{$i$}}
\put(101.68,5.33){\makebox(0,0)[cc]{$j$}}
\put(79.35,13.33){\vector(0,1){0.2}}
\emline{79.35}{9.66}{33}{79.35}{13.33}{34}
\put(85.68,5.33){\makebox(0,0)[cc]{$s$}}
\put(95.35,5.33){\makebox(0,0)[cc]{$t$}}
\put(57.68,17.66){\makebox(0,0)[lc]{,\quad $\etat_{istj}=$}}
\put(77.01,17.66){\framebox(11.67,7.00)[cc]{$\pda$}}
\put(92.68,17.66){\framebox(11.67,7.00)[cc]{$\pch$}}
\put(85.34,14.33){\vector(0,-1){0.2}}
\emline{85.34}{17.66}{35}{85.34}{14.33}{36}
\put(85.35,13.33){\vector(0,1){0.2}}
\emline{85.35}{9.66}{37}{85.35}{13.33}{38}
\put(95.00,17.33){\vector(0,1){0.2}}
\emline{95.00}{13.00}{39}{95.00}{17.33}{40}
\put(48.33,9.66){\vector(0,-1){0.2}}
\emline{48.33}{13.00}{41}{48.33}{9.66}{42}
\put(95.00,9.66){\vector(0,-1){0.2}}
\emline{95.00}{13.00}{43}{95.00}{9.66}{44}
\put(101.34,9.66){\vector(0,-1){0.2}}
\emline{101.34}{13.00}{45}{101.34}{9.66}{46}
\put(108.67,15.00){\makebox(0,0)[cc]{.}}
\end{picture}
\\
Figure\,8: The elements $\xi^\tau_{istj}$ and $\etat_{istj}$.
\end{center}
By \rf[e-vwst] one has
\begin{center}
\special{em:linewidth 0.4pt}
\unitlength 1.00mm
\linethickness{0.4pt}
\begin{picture}(58.00,35.33)
\put(22.66,26.00){\makebox(0,0)[rc]{$\etat_{st}\rac u^i_j=$}}
\put(26.66,26.00){\framebox(11.67,7.00)[cc]{$\pda$}}
\put(42.33,26.00){\framebox(11.67,7.00)[cc]{$\pch$}}
\put(44.66,35.33){\vector(0,1){0.2}}
\emline{44.66}{33.33}{1}{44.66}{35.33}{2}
\put(50.99,35.33){\vector(0,1){0.2}}
\emline{50.99}{33.33}{3}{50.99}{35.33}{4}
\put(29.33,33.33){\vector(0,-1){0.2}}
\emline{29.33}{35.33}{5}{29.33}{33.33}{6}
\put(35.33,33.33){\vector(0,-1){0.2}}
\emline{35.33}{35.00}{7}{35.33}{33.33}{8}
\put(28.99,22.67){\vector(0,-1){0.2}}
\emline{28.99}{26.00}{9}{28.99}{22.67}{10}
\put(50.99,25.67){\vector(0,1){0.2}}
\emline{50.99}{21.34}{11}{50.99}{25.67}{12}
\put(28.99,13.67){\makebox(0,0)[cc]{$s$}}
\put(56.66,13.67){\makebox(0,0)[cc]{$j$}}
\put(51.00,14.00){\makebox(0,0)[cc]{$t$}}
\put(56.66,17.67){\vector(2,-1){0.2}}
\emline{52.00}{19.67}{13}{56.66}{17.67}{14}
\put(29.00,19.00){\vector(0,1){0.2}}
\emline{29.00}{17.34}{15}{29.00}{19.00}{16}
\put(23.00,13.67){\makebox(0,0)[cc]{$i$}}
\emline{23.00}{17.34}{17}{25.15}{18.83}{18}
\emline{25.15}{18.83}{19}{27.31}{20.13}{20}
\emline{27.31}{20.13}{21}{29.47}{21.24}{22}
\emline{29.47}{21.24}{23}{31.65}{22.17}{24}
\emline{31.65}{22.17}{25}{33.83}{22.91}{26}
\emline{33.83}{22.91}{27}{36.03}{23.45}{28}
\emline{36.03}{23.45}{29}{38.23}{23.81}{30}
\emline{38.23}{23.81}{31}{41.33}{24.00}{32}
\put(43.66,25.00){\vector(1,4){0.2}}
\emline{38.33}{22.00}{33}{39.99}{21.25}{34}
\emline{39.99}{21.25}{35}{41.36}{21.22}{36}
\emline{41.36}{21.22}{37}{42.45}{21.92}{38}
\emline{42.45}{21.92}{39}{43.66}{25.00}{40}
\emline{35.33}{26.00}{41}{36.00}{24.00}{42}
\emline{45.00}{23.34}{43}{47.33}{22.77}{44}
\emline{47.33}{22.77}{45}{50.00}{21.00}{46}
\put(58.00,26.67){\makebox(0,0)[cc]{.}}
\put(51.00,18.00){\vector(0,-1){0.2}}
\emline{51.00}{21.34}{47}{51.00}{18.00}{48}
\end{picture}
\\
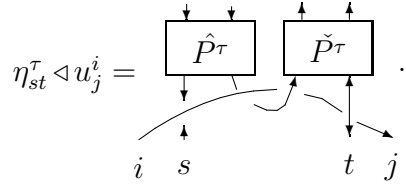
Figure\,9: The right adjoint action of $u^i_j$  on $\etat_{st}$.
\end{center}
Replacing one crossing $\Rda{}$ by $\Rdam{}+\QM I-\QM K$ similarly to
the graphical calculations in the first part of the proof one can show that
\begin{align*}
\zeta^\tau_{sijt}:=\etat_{uv}\rac
u^a_b\,\Rch[au]{si}\Rda[vb]{jt}-\delta_{ij}\etat_{st}=\QM \xi^\tau_{sijt}-\QM\rtm\etat_{sijt}.
\end{align*}
Since $q$ is not a root of unity, $T=\QM(\Rdam{}-\rtm I)$ is invertible with
inverse $T^{-1}=\QM^{-1}\bigl((\qm-\rtm)^{-1}\pda[+]-(q+\rtm)^{-1}\pda[-]\bigr)$.
Therefore
$$
B^y_i(T^{-1})^{yj}_{zl}C^k_z\zeta^\tau_{sijt}=\etat_{sklt}
$$
belongs to $\Lam[\tau]_\ll$. Finally we have
$B^m_vB^a_yC^b_zC^l_t\etat_{vyzt}\Rre[ab]{nk}=\theta^\tau_{mnkl}$ which
completes the proof.
\end{proof}
Now we are ready to complete the proof of the theorem. By
Lemma\,\ref{l-biinv}\,(iii) both $\etat[+] $ and $\etat[-]$ belong to
$\SS(\RR)$ (see Section\,\ref{uni}) and $\dim\gdui[2]\ge1$. Hence
$\eta^0\not\in\SS(\RR)_\ii$. Since $\theta\ott\theta\equiv\rtm\etat[0]\mod
\SS(\RR)$   by
\rf[e-trans] and  $a\etat[0]=\etat[0] a$, $a\in\AA$, by
\rf[e-tn] we get
\begin{align}\label{e-central}
a\theta{\land}\theta=\theta{\land}\theta a,\quad a\in\AA.
\end{align}
We prove \rf[e-inner] by induction over  the degree $n$ of $\rho\in\gdu[n]$. For
$n=0$ it is true by the definition of the FODC. Suppose it is true for
$n-1$. Since there exist $\alpha_i\in\gdu[n-1]$ and $b_i\in\AA$ such that
$\rho=\alpha_i \dd b_i$, we obtain by induction assumption and by \rf[e-central]
\begin{align*}
\dd\rho=\dd \alpha_i \dd b_i&=
\theta \alpha_i\dd b_i-(-1)^{n-1}\alpha_i \theta(\theta b_i-b_i\theta)
\\
&=\theta\rho-(-1)^{n}(\alpha_i\theta b_i\theta-\alpha_i b_i\theta^2)
\\
&=\theta\rho-(-1)^n\alpha_i\dd b_i\theta
\\
&=\theta\rho-(-1)^n\rho\theta.
\end{align*}
Using $\dd^2\rho=0$ and \rf[e-inner] twice gives  $\theta^2\rho=\rho\theta^2$,
and $\theta^2$ is central in $\gdu[2]$. This completes the proof of (ii).
\\
By  Lemma\,\ref{l-biinv}\,(iii) and Lemma\,\ref{l-generator} 
$$
\SS(\RR)\supseteq \SS(\RR)_0\rac\AA =\Lam[+]_\ll\oplus\Lam[-]_\ll.
$$
Since $\SS(\RR)\subseteq\ker A_2$ by universality of $\gdu$ and $\eta^0\in\ker
A_2$, we conclude with Lemma\,\ref{l-sub}\,(ii) that the above inclusion is not strict, $\SS(\RR)=\Lam[+]_\ll{\oplus}\Lam[-]_\ll$ and
$$
(\ker A_2)_\ll=\SS(\RR)\oplus \Lam[0]_\ll.
$$
Since both $\gdul$ and $\gdsl$ are quadratic algebras,
$\gdul{/}(\etat[0])\cong\gten_\ll{/}(\SS(\RR)\oplus\langle\etat[0]\rangle)\cong
\gten_\ll{/}(\ker A_2)_\ll\cong\gdsl$. Since both $\gdu{/}(\etat[0])$ and $\gds$ are {\em free}
left $\AA$-modules it follows $\gdu{/}(\etat[0])\cong\gds$. Noting that
$\theta^2=\rtm\etat[0]$ in $\gdu$ completes the proof of the theorem.

\section{The bi-invariant 2-form of the universal differential calculus}\label{uni}

In this section we will complete the proof of Lemma\,\ref{l-biinv}\,(iii) and show that both bi-invariant elements $\eta^+$ and
$\eta^-$ belong to $\SS(\RR)$.  We give  different proofs for the cases $\Gamm_+$ and $\Gamm_-$.
The first proof for $\Gamm_+$ is self-contained and  much easier than the
second one. In the later one  we take  results from \cite{a-SchSch1} and make
use of  a computer algebra program to simplify long terms. For $q$ {\em
  transcendental} however the first proof works for $\Gamm_-$ as well. 
\\
We recall some identities  which are easily proved using \rf[e-om] and the
Leibniz rule. Equations \rf[e-omab],
\cite[formula (14.3)]{b-KS},  and \rf[e-ssab], \cite[Lemma\,14.15]{b-KS},
 are valid for
arbitrary {\em left-covariant} FODC while \rf[e-omrac] and \rf[e-ssa] in addition require $\dd
a=\theta a-a\theta$. For $a,b\in\AA$ and $p\in\RR$ we have
\begin{align}\label{e-omrac}
\om(a)&=\theta\rac a+\ve(a)\theta,\qquad\theta\rac p=0,
\\
\om(ab)&=\om(a)\rac b+\ve(a)\om(b),\label{e-omab}
\\
\begin{split}
\SS(a)&=(\theta\ott\theta)\rac a-\theta\ot(\theta\rac a)-(\theta\rac a)\ot 
\theta+\ve(a)\theta\ot\theta,
\\
 \SS(p)&=(\theta\ott\theta)\rac p,\label{e-ssa}
\end{split}
\\
\begin{split}
\SS(\plus{a}b)&=\SS(a)\rac b+\om(a)\rac b_{(1)}\ot\om(b_{(2)})+\om(b_{(1)})\ot 
(\om(a)\rac b_{(2)}),\label{e-ssab}
\\
\SS(pb)&=\SS(p)\rac b.
\end{split}
\end{align}

We abbreviate $\RM=\rt-\rtm$. In what follows we do {\em not} sum over signs
$\tau$ and $\nu$.

{\em Part~1. $\Gamm=\Gamm_+$.}  First
we show $\Q:=\plus{U}{\cdot}\plus{U}-\QM\RM\plus{U}\in\RR$. 
By \rf[e-thetau] and \rf[e-trq] we obtain
$
\theta\rac U=\trq(\Rda[2]{})^m_n\theta_{mn}=(\QM\RM+\xx)\theta.
$ Using \rf[e-omrac], \rf[e-omab],  and $\ve(U)=\xx$ we have 
$\om(\Q)=(\theta\rac\plus{U})\rac\plus{U}-\QM\RM\theta\rac\plus{U}=0$. In addition
$\ve(\Q)=0$; hence $Q\in\RR$. Next we compute $\SS(\Q)$. 
Since $\Q\in\RR$, by \rf[e-ssa]   we have
$\SS(\Q)=(\theta\ott\theta)\rac\Q$. Using \rf[e-trans] and \rf[e-etatuu] we get
\begin{align}
\SS(\Q)&=(\theta\ott\theta)\rac\Q=(q\etat[+]-\qm\etat[-])\rac(\plus{U}\plus{U}-\QM\RM\plus{U})\nn
\\
&=q\alpt[+]^2\etat[+]-\qm\alpt[-]^2\etat[-]-\QM\RM(q\alpt[+]\etat[+]-\qm\alpt[-]\etat[-])\nn
\\
&=q\alpt[+](\alpt[+]-\QM\RM)\etat[+]-\qm\alpt[-](\alpt[-]-\QM\RM)\etat[-].\label{e-1}
\end{align}
Since $\RR$ is a right ideal $\Q\plus{U}\in\RR$. By  \rf[e-ssab] and \rf[e-1] 
\begin{align}\label{e-2}
\SS(\Q\plus{U})=q\alpt[+]^2(\alpt[+]-\QM\RM)\etat[+]-\qm\alpt[-]^2(\alpt[-]-\QM\RM)\etat[-].
\end{align}
Solving this linear system \rf[e-1] and \rf[e-2] in $\etat[+]$ and $\etat[-]$
we have to consider its  coefficient determinant 
\begin{align*}
\det&=\alpt[+]\alpt[-](\alpt[+]-\alpt[-])(\alpt[+]-\QM\RM)(\alpt[-]-\QM\RM)
\\
&=(\rt+\rtm)\QM^6\qtwo^3(q\rt-\qm\rtm)(\qm\rt-q\rtm)(q^2\rt-q^{-2}\rtm)(q^{-2}\rt-q^2\rtm).
\end{align*}
Since $q$ is not a root of unity, $\det\ne0$. Hence both $\etat[+]$ and $\etat[-]$
 belong to $\SS(\RR)$.

{\em Part~2. $\Gamm=\Gamm_-$}. We denote the critical value by $\cv$, 
$\cv=\QM q^2\rt(2\xx+\QM\RM)=(q^4+1)\rt^2+2q(q^2-1)\rt-(q^4+1)$.
We recall some of the  defining constants for $\Gamm_-$ from \cite[p.\,656]{a-SchSch1} 
\begin{align*}
\mu^+ &=\frac{\RM(-q^2\rt+q^{-2}\rtm-\QM)}{\QM\RM+2\xx},
\\
\mu^- &=\frac{\RM(-q^{-2}\rt+q^2\rtm-\QM)}{\QM\RM+2\xx}.
\end{align*}
The idempotents $\pda[\nu]$, $\nu\in\{+,-\}$, and their $q$-traces are as follows
\begin{align}\label{e-pro}
\pda[\nu]&=(\lamn+\lamn^{-1})^{-1}(\lamn^{-1} I+\Rda{}+\QM(1-\rt\lamn)^{-1}K),
\\
\trq(\pda[\nu])&=\frac{\RM(\lamn^2\rt-\lamn^{-1})}{\QM(\lamn+\lamn^{-1})(\lamn\rt-1)}I=:\tnn
I.\label{e-tpro}
\end{align}
There are two $\Adr$-invariant quadratic elements in $\AA$, namely 
$V_\nu =D^b_aD^j_i(\pda[\nu])^{ai}_{yz}u^y_bu^z_j$. One has
$\ve(V_\nu)=\xx\tnn$.
Note that $\plus{W_\nu}\in\RR$, where $W_\nu=V_\nu-\mu^\nu U$, \cite[p.\,656,
 eq.\,(3)]{a-SchSch1}. 
Suppose $a\in\AA$ is $\Adr$-invariant and $\rho\in\Lam_\ii$, where $\Lam$ is
a bicovariant bimodule. Then one has $\rho\rac a\in\Lam_\ii$. Namely, 
$\dr(\rho\rac a)=(Sa_{(2)}\ot Sa_{(1)})(\rho\ot1)(a_{(3)}\ot a_{(4)})=\rho\rac a_{(2)}\ot Sa_{(1)}a_{(3)}=\rho\rac a\ot 1.$
Applying this fact to $\Lam[\tau]$ and  $V_\nu$, and noting
that $\etat$, $\tau\in\{+,-\}$, is the only bi-invariant element of
$\Lam[\tau]$ (up to scalars), there exist complex numbers  $\ctn$, $\tau,\,\nu\in\{+,-\}$, defined by
\begin{align*}
\etat\rac\plus{V_\nu}=\ctn\etat.
\end{align*}
We shall determine these constants.
By the definition of $V_\nu$ and \rf[e-etatu]
\begin{align*}
\etat\rac V_\nu&=\etat\rac u^i_ju^s_t(\pda[\nu])^{jt}_{ab}D^a_iD^b_s
\\
&=\bigl(\delta_{ij}\etat+ \QM(\lam_\tau^2+1)\xi^\tau_{ij}-
\rtm\QM(1+\lam_\tau^{-2})\etat_{ij}\bigr)\rac u^s_t(\pda[\nu])^{jt}_{ab}D^a_iD^b_s.
\end{align*}
We carry out the calculations for the first term. By \rf[e-etatu] and \rf[e-tpro] we have
\begin{align*}
\delta_{ij}\etat&\rac(u^s_t(\pda[\nu])^{jt}_{ab}D^a_iD^b_s)=\tnn\etat\rac U=\tnn (\xx+\alpt)\etat.
\end{align*}
Using  graphical calculations we obtain for the other two terms
\begin{align*}
\xi^\tau_{ij}&\rac(u^s_t(\pda[\nu])^{jt}_{ab}D^a_iD^b_s)=(\QM\delta_{\tau,\nu}\rt^2\lamt^2+\rt\lamn^2\tnn-\QM\lamn\lamt^{-1}e_{\tau\nu})\etat,
\\
\etat_{ij} & \rac(u^s_t(\pda[\nu])^{jt}_{ab}D^a_iD^b_s)=(\QM\rt\lamn^{-1}\lamt e_{\tau\nu}+\lamn^{-2}\tnn -\QM\rtm\lamt^{-2}\delta_{\tau,\nu})\etat,
\end{align*}
where $e_{\tau\nu}=(\lamn+\lamn^{-1})^{-1}(\lamt^{-1}+\QM(1-\rt\lamn)^{-1})$ is 
obtained from \rf[e-pro] and  the picture
\begin{center}
\special{em:linewidth 0.4pt}
\unitlength 1.00mm
\linethickness{0.4pt}
\begin{picture}(67.00,26.99)
\put(23.99,17.66){\framebox(11.67,7.00)[cc]{$\pda$}}
\put(39.66,17.66){\framebox(11.67,7.00)[cc]{$\pch$}}
\put(41.99,26.99){\vector(0,1){0.2}}
\emline{41.99}{24.99}{1}{41.99}{26.99}{2}
\put(48.32,26.99){\vector(0,1){0.2}}
\emline{48.32}{24.99}{3}{48.32}{26.99}{4}
\put(26.66,24.99){\vector(0,-1){0.2}}
\emline{26.66}{26.99}{5}{26.66}{24.99}{6}
\put(32.66,24.99){\vector(0,-1){0.2}}
\emline{32.66}{26.66}{7}{32.66}{24.99}{8}
\put(32.00,7.66){\framebox(11.67,7.00)[cc]{$\acute{P^\nu}$}}
\put(34.33,15.00){\vector(2,-3){0.2}}
\emline{32.67}{17.67}{9}{34.33}{15.00}{10}
\put(41.67,17.33){\vector(1,2){0.2}}
\emline{40.67}{15.00}{11}{41.67}{17.33}{12}
\put(34.33,7.33){\vector(2,3){0.2}}
\emline{26.67}{17.33}{13}{26.74}{14.46}{14}
\emline{26.74}{14.46}{15}{26.88}{11.90}{16}
\emline{26.88}{11.90}{17}{27.07}{9.65}{18}
\emline{27.07}{9.65}{19}{27.33}{7.72}{20}
\emline{27.33}{7.72}{21}{27.64}{6.09}{22}
\emline{27.64}{6.09}{23}{28.02}{4.78}{24}
\emline{28.02}{4.78}{25}{28.46}{3.78}{26}
\emline{28.46}{3.78}{27}{28.96}{3.09}{28}
\emline{28.96}{3.09}{29}{29.51}{2.71}{30}
\emline{29.51}{2.71}{31}{30.13}{2.65}{32}
\emline{30.13}{2.65}{33}{30.81}{2.89}{34}
\emline{30.81}{2.89}{35}{31.55}{3.45}{36}
\emline{31.55}{3.45}{37}{32.35}{4.32}{38}
\emline{32.35}{4.32}{39}{33.22}{5.50}{40}
\emline{33.22}{5.50}{41}{34.33}{7.33}{42}
\put(48.33,17.33){\vector(0,1){0.2}}
\emline{41.00}{7.33}{43}{42.00}{5.77}{44}
\emline{42.00}{5.77}{45}{42.94}{4.50}{46}
\emline{42.94}{4.50}{47}{43.80}{3.52}{48}
\emline{43.80}{3.52}{49}{44.58}{2.83}{50}
\emline{44.58}{2.83}{51}{45.30}{2.44}{52}
\emline{45.30}{2.44}{53}{45.94}{2.33}{54}
\emline{45.94}{2.33}{55}{46.50}{2.52}{56}
\emline{46.50}{2.52}{57}{47.00}{3.00}{58}
\emline{47.00}{3.00}{59}{47.42}{3.77}{60}
\emline{47.42}{3.77}{61}{47.77}{4.83}{62}
\emline{47.77}{4.83}{63}{48.04}{6.19}{64}
\emline{48.04}{6.19}{65}{48.25}{7.83}{66}
\emline{48.25}{7.83}{67}{48.38}{9.77}{68}
\emline{48.38}{9.77}{69}{48.43}{12.00}{70}
\emline{48.43}{12.00}{71}{48.42}{14.52}{72}
\emline{48.42}{14.52}{73}{48.33}{17.33}{74}
\put(67.00,14.67){\makebox(0,0)[cc]{$=e_{\tau\nu}\etat.$}}
\end{picture}
\\
Figure\,10: The definition of $e_{\tau\nu}$.
\end{center}
Finally we obtain
\begin{align*}
\begin{split}
\ctn[++]&=\cv^{-1}(q^2+1)(q^2 \rt^2-1)\rt^{-2}q^{-6}((q^{12}+q^4) \rt^4  +(2 q^{11} -2 q^9 +2 q^5 -2
)\rt^3+
\\
&\qquad+(-q^{12} +q^8 -4 q^6 +q^4 -1)\rt^2+(-2 q^9 +2 q^7 -2 q^3 +2 q)\rt+q^8+1)
\end{split}
\\
\ctn[+-]&=\cv^{-1}\qtwo(2(q-\qm\rt^{-2})(\rt+q^3)(\rt-q)((q^2 -1 +q^{-2})\rt^2+\QM\rt -q^2+1-q^{-2})
\\
\ctn[-+]&=\cv^{-1}\qtwo(2(\qm-q\rt^{-2})(q \rt+1)(q^3 \rt-1)((q^2 -1 +q^{-2})\rt^2+\QM\rt-q^2+1-q^{-2})
\\
\begin{split}
\ctn[--]&=\cv^{-1}(q^2+1)(\rt^2-q^2)\rt^{-2}q^{-6}((q^8 +1)\rt^4  
+(2 q^9 -2 q^7 +2 q^3 -2 q )\rt^3+
\\
&\qquad+(-q^{12}+q^8 -4 q^6 +q^4 -1)\rt^2+(-2 q^{11} +2 q^9 -2 q^5 +2 q^3)\rt+q^{12}+q^4)
\end{split}
\end{align*}
Now we are able to compute four bi-invariant elements of $\SS(\RR)$.
Using \rf[e-trans]  and \rf[e-etatuu] we have
\begin{align}
\SS(W_\nu)&=(\theta\ott\theta)\rac (\plus{V_\nu}-\mu^\nu \plus{U})=q(\ctn[+\nu]-\mu^\nu\alpt[+])\etat[+]-\qm(\ctn[-\nu]-\mu^\nu\alpt[-])\etat[-].\label{e-3}
\end{align}
Similarly to Part~1 we get
\begin{align}\label{e-4}
\SS(W_\nu\plus{U})=q\alpt[+](\ctn[+\nu]-\mu^\nu\alpt[+])\etat[+]-\qm\alpt[-](\ctn[-\nu]-\mu^\nu\alpt[-])\etat[-].
\end{align}
Consider the ${4{\times} 2}$--coefficient matrix  $T=(T_{ij})$  for the
linear system of equations \rf[e-3] and
\rf[e-4], $\nu\in\{+,-\}$. The two columns are
\begin{align*}
(T_{i1})&=q\bigl(\ctn[++]-\mu^+\alpt[+],\,\,\ctn[+-]-\mu^-\alpt[+],\,\,\alpt[+](\ctn[++]-\mu^+\alpt[+]),\,\,\alpt[+](\ctn[+-]-\mu^-\alpt[+])\bigr)^\fettt,
\\
(T_{i2})&=-\qm\bigl(\ctn[-+]-\mu^+\alpt[-],\,\,\ctn[--]-\mu^-\alpt[-],\,\,\alpt[-](\ctn[-+]-\mu^+\alpt[-]),\,\,\alpt[-](\ctn[--]-\mu^-\alpt[-])\bigr)^\fettt.
\end{align*}
We distinguish three cases: $(T_{i1})=0$,  $(T_{i2})=0$, and no column
vanishes, respectively.

{\em Case~1.} $\ctn[++]-\mu^+\alpt[+]=\ctn[+-]-\mu^-\alpt[+]=0$.
We obtain the following two equations
\begin{align*}
\begin{split}
0&= (q^2+1)(q\rt-1)(q\rt+1)(q^{12}\rt^4+q^4\rt^4+2q^{11}\rt^3-2q^9\rt^3+2q^5\rt^3-2q^3\rt^3-
\\
&\quad -q^{12}\rt^2+q^8\rt^2-4q^6\rt^2+q^4\rt^2-\rt^2-2q^9\rt+2q^7\rt-2q^3\rt+2q\rt+q^8+1),
\end{split}
\\
0&=2(q\rt-1)(q\rt+1)(\rt+q^3)(\rt-q)\times
\\
&\quad\quad \times(q^4\rt^2-q^2\rt^2+\rt^2+q^3\rt-q\rt-q^4+q^2-1)(q^2+1).
\end{align*}
Since $q$ is not a root of unity we have
\begin{align*}
\begin{split}
d_1&= q^{12}\rt^4+q^4\rt^4+2q^{11}\rt^3-2q^9\rt^3+2q^5\rt^3-2q^3\rt^3-
\\
&\quad -q^{12}\rt^2+q^8\rt^2-4q^6\rt^2+q^4\rt^2-\rt^2-2q^9\rt+2q^7\rt-2q^3\rt+2q\rt+q^8+1=0,
\end{split}
\\
d_2&=q^4\rt^2-q^2\rt^2+\rt^2+q^3\rt-q\rt-q^4+q^2-1=0.
\end{align*}
Using the Euclidean algorithm we eliminate powers of $\rt$.  
We end up with polynomials 
\begin{align*}
a=& (q^4+1)(q^6-q^2+1)(q^8+1)q^4\rt^3+
\\
& +(q^4-q^2+1)(3q^{14}-2q^{10}+2q^8-q^6+2q^2-2)q^3\rt^2+
\\
& +(q^{20}-2q^{18}+2q^{14}-4q^{12}-q^{10}+5q^8-6q^6+q^4+q^2-1)\rt+
\\
&+(q^4-q^2+1)(q^{16}-4q^{14}+2q^{12}+2q^{10}-5q^8+2q^6+2q^4-4q^2+2)q
\end{align*}
 and 
\begin{align*}
b&=(-q^4-1)(q^6-q^2+1)(q^4-q^2+1)^4(q-1)^6(q+1)^6q\rt+
\\
&\quad\quad+ (-2q^{12}+3q^{10}-3q^8+q^6+q^4-2q^2+1)(q^4-q^2+1)^3(q-1)^6(q+1)^6q^2 
\end{align*}
such that $a d_2+b d_1=-(q^6+q^3+1)(q^6-q^3+1)(q-1)^{6}(q+1)^{6}q$ (there is no $\rt$ left). Since $d_1=d_2=0$, $q$ is a root of unity which contradicts our assumption. 
Hence Case~1 is impossible.

{\em Case~2.} $\ctn[-+]-\mu^+\alpt[-]=\ctn[--]-\mu^-\alpt[-]=0$.
Similarly to Case~1 we have
\begin{align*}
d_3&=q^4\rt^2-q^2\rt^2+\rt^2+q^3\rt-q\rt-q^4+q^2-1=0,
\\
\begin{split}
d_4&=q^8\rt^4+\rt^4+2q^9\rt^3-2q^7\rt^3+2q^3\rt^3-2q\rt^3-
\\
 &\quad -q^{12}\rt^2+q^8\rt^2-4q^6\rt^2+q^4\rt^2-\rt^2-2q^{11}\rt+2q^9\rt-2q^5\rt+2q^3\rt+q^{12}+q^4=0.
\end{split}
 \end{align*}
Again there exist polynomials $a$ and $b$ in $q$ and $\rt$ such that $ad_3+bd_4={(q^6+q^3+1)}{(q^6-q^3+1)}{(q^4-q^2+1)}{(q-1)^{12}}{(q+1)^{12}}q^2$. 
This  contradicts our assumption that $q$ is not a root of unity. Hence,
the only possibility is 

{\em Case~3.} We will show, that $T$ has rank $2$. Suppose to the contrary that $T$ has at least rank $1$. Then the ${2{\times}2}$-matrices  formed by the first and third rows, respectively from the second and fourth rows, both have zero determinant. Since $\alpha_+-\alpha_-\ne0$ this is equivalent to
 $(\ctn[++]-\mu^+\alpt[+])(\ctn[-+]-\mu^+\alpt[-])=0$ and $(\ctn[+-]-\mu^-\alpt[+])(\ctn[--]-\mu^-\alpt[-])=0$. Since moreover the matrix
 formed by  the first two rows has zero determinant, we conclude
 $\ctn[++]-\mu^+\alpt[+]=\ctn[+-]-\mu^-\alpt[+]=0$ or $\ctn[-+]-\mu^+\alpt[-]=\ctn[--]-\mu^-\alpt[-]=0$. But this is impossible by
 Case~1 and Case~2. Hence $T$ has rank $2$; both $\eta^+$ and $\eta^-$ belong
 to $\SS(\RR)$. 

{\em Acknowledgements.} The author is greatly indebted to I.\,Heckenberger for
many stimulating discussions and helpful comments.

\providecommand{\bysame}{\leavevmode\hbox to3em{\hrulefill}\thinspace}

\end{document}